\documentclass[11pt,dvipsnames,reqno]{amsart}

\usepackage{times}
\usepackage[all]{xy}
\usepackage{amssymb,amsmath,amsthm}
\usepackage[dvipdfmx]{graphicx}

\newcommand{\calO}{\mathcal{O}}
\newcommand{\calP}{\mathcal{P}}
\newcommand{\calN}{\mathcal{N}}

\newcommand{\calF}{\mathcal{F}}

\newcommand{\bbP}{\mathbb{P}}
\newcommand{\bbZ}{\mathbb{Z}}
\newcommand{\bbQ}{\mathbb{Q}}

\newcommand{\bbR}{\mathbb{R}}
\newcommand{\bbC}{\mathbb{C}}

\newcommand{\bfF}{\mathbf{F}}
\newcommand{\bbF}{\mathbb{F}}

\newcommand{\frakG}{\mathfrak{G}}

\newcommand{\Pic}{\mathrm{Pic}}
\usepackage{graphicx}
\newcommand{\half}{\frac{1}{2}}
\newcommand{\Aut}{\mathrm{Aut}}
\newcommand{\Hom}{\mathrm{Hom}}
\newcommand{\rank}{\mathrm{rank}}

\newcommand{\Sing}{\mathrm{Sing}}

\usepackage[
           colorlinks=true,
            urlcolor=blue,       
            filecolor=blue,     
            linkcolor=blue,       
            citecolor=blue,         
            pdftitle= {Title},
            pdfauthor={Author},
            pdfsubject={Subject},
            pdfkeywords={Key}
            pagebackref,
            pdfpagemode=None,
            bookmarksopen=true]
            {hyperref}

\newcommand{\frakS}{\mathfrak{S}}

\newtheorem{theorem}{Theorem}[section]
 \newtheorem{proposition}[theorem]{Proposition}
 \newtheorem{lemma}[theorem]{Lemma}
 \newtheorem{corollary}[theorem]{Corollary}

 \theoremstyle{definition}
 
 \newtheorem{defn}[theorem]{Definition}
 \newtheorem{remark}[theorem]{Remark}

\newcommand{\Kum}{\mathrm{Kum}}

\newcommand{\Bl}{\operatorname{Bl}}

\newcommand{\la}{\langle}
\newcommand{\ra}{\rangle}

\newcommand{\beq}{\begin{equation}}
\newcommand{\eeq}{\end{equation}}

\def\MARU#1{\textcircled{\scriptsize#1}}

\begin{document}
\title [Desmic quartic surfaces in arbitrary characteristic]{Desmic quartic surfaces in arbitrary characteristic}

\author{Igor Dolgachev}
\address{Department of Mathematics, University of Michigan, 525 East University Avenue, Ann Arbor, MI 48109-1109 USA}
\email{idolga@umich.edu}

\author{Shigeyuki Kond\=o}
\address{Graduate School of Mathematics, Nagoya University, Nagoya, 464-8602 Japan}
\email{kondo@math.nagoya-u.ac.jp}
\thanks{Research of the second author is partially supported by JSPS Grant-in-Aid for Scientific Research (A) No.20H00112, (B) No.25K00906.}

\maketitle

\begin{center}
{\footnotesize \it Dedicated to Tetsuji Shioda on the occation of his 85th birthday}
\end{center}

\begin{abstract}
A desmic quartic surface is a birational model of the Kummer surface of the self-product of an elliptic curve.
We recall the classical geometry of these surfaces and study their analogs in arbitrary characteristic.
Moreover, we discuss the cubic line complex $\frakG$ associated with the desmic tetrahedra introduced by G. Humbert.
We prove that $\frakG$ is a rational Fano threefold with $34$ nodes. 
The number $34$ is the maximum number
of nodes on a Fano threefold of degree 6 in $\bbP^5$,
and the group of projective automorphisms is isomorphic
to $\frakS_4\wr 2 = (\frakS_4\times \frakS_4)\rtimes 2$.
\end{abstract}

\section{Introduction} 
Two tetrahedra $T_1$ and $T_2$ in a real or complex three-dimensional projective space are called \emph{desmic} if each 
edge of one intersects a pair of opposite edges of another \cite[Section 3.19]{Lord}. Equivalently, $T_1$ and $T_2$ are desmic 
if the edges of a face of one 
are the diagonals of the complete quadrilateral formed by the intersection of the face with another tetrahedron. 
Two desmic tetrahedra are perspective from each vertex of another tetrahedron 
$T_3$, and each pair of three tetrahedra $T_1, T_2$, and $T_3$ are desmic.  The three tetrahedra form a 
\emph{desmic system} of tetrahedra. Three tetrahedra form a desmic system if and only if they belong to the same pencil 
of quartic surfaces, which is called a \emph{desmic pencil}. All its members, different from the three tetrahedra, 
are irreducible quartic surfaces 
with 12 nodes lying by pairs in the edges of the tetrahedra that intersect a pair of opposite edges of a desmic 
tetrahedron. An irreducible member of a desmic pencil is called a \emph{desmic quartic surface}.

 Desmic tetrahedra were first introduced and studied 
 by Cyparissos Stephanos \cite{Stephanos} in 1879.\footnote{$\delta\epsilon\sigma\mu$\'o$\varsigma$ = bond} In 1891, Georges Humbert showed that an irreducible member of the pencil of quartic surfaces 
containing three desmic tetrahedra is birationally isomorphic to the Kummer surface of the self-product of an elliptic 
curve \cite{Humbert}. He also showed that a desmic quartic surface is isomorphic to Cremona's quartic surface, 
the locus of nodes of quadric surfaces that cut out the union of three conics in a smooth cubic surface 
whose spanning planes contain three coplanar lines on the surface. Other beautiful geometric properties of desmic quartic surfaces were later studied in two papers 
by Mathews \cite{Mathews1}, \cite{Mathews2}, and in Jessop's book on quartic surfaces \cite{Jessop}.

A desmic quartic surface is singular. It contains 12 ordinary nodes, two on each edge of a desmic tetrahedron. 
It also contains 16 lines, which, together with the nodes, form an abstract configuration 
$(12_4,16_3)$ (see the definition in Section 3) isomorphic 
to the Reye or the Hesse-Salmon configurations \cite{DolgachevAbstract}. 
If the characteristic $p\ne 2$, any surface containing such a configuration is 
birationally isomorphic to the Kummer surface $\Kum(E\times E)$ of the self-product of 
an elliptic curve. 
 
In this paper, we first recall the classical results about desmic quartic surfaces 
and study the cubic line complex associated with the pencil first introduced by Humbert. 
We find its explicit equation and show that it has 34 ordinary nodes and 24 planes. We associate a cubic $5$-fold with this cubic line complex, which has the maximum number of nodes, 35.  It follows that
the number 34 of nodes is the maximum number 
of terminal singularities on a Fano threefold of degree 6 in $\bbP^5$.  
This cubic 5-fold is projectively isomorphic to the Segre cubic 5-fold with
35 nodes (see Theorem \ref{maximum}). 
We also show the rationality of the cubic line complex.
The cubic complex has the group of projective automorphisms isomorphic to 
the group $\frakS_4\wr 2$, and can be 
viewed as the analog of the Segre cubic primal with 10 nodes and 15 planes, and the group of projective 
symmetries isomorphic to 
$\frakS_6$. The group $\frakS_4\wr 2$ is closely related to the Weyl group of the root system of 
type $F_4$, known to be isomorphic to the group of symmetries of the $24$-cell. 
We show that the quotient of the Weyl group  
$W(F_4)$ by its center is isomorphic to the group of projective symmetries of the desmic pencil. We also 
give a geometric interpretation of the subgroups $W(D_4)$ and $W(B_4)$ of $W(F_4)$.

Then, we study the analogs of desmic quartic surfaces in characteristic two. One of them could be the Kummer surface of the self-product 
 of an elliptic curve.   
 Instead of $12$ nodes, the surface has six rational double points of type $A_3$ in the case of an ordinary elliptic curve. We don't know the case of the supersingular elliptic curve. 
 Another analog is Cremona's quartic surface, which can be 
 defined 
 in any characteristic, but acquires an additional rational double point if the characteristic is equal to $2$.
In this case, we show that the surface is birationally isomorphic to a supersingular K3 surface with Artin invariant 
 $\sigma\le 2$. The one-dimensional family of such surfaces contains a supersingular surface with Artin invariant one, 
 previously studied in the paper \cite{DK}, and we show how to find the configuration of $(12_4,16_3)$ among 
 the set of smooth rational curves on this surface explicitly. 
 
 \vskip4pt
The authors thank Ivan Cheltsov and Shigeru Mukai for a fruitful discussion about the geometry of the cubic 
line complex. We are also grateful to the referee for many helpful comments that substantially improved the paper's exposition.

\section{Desmic tetrahedra}
Until Section 6, we assume that the ground field $\Bbbk$ is the field $\bbC$ of complex numbers, although many of the 
discussions are valid over an algebraically closed field of characteristic different from $2$. 
A \emph{tetrahedron} means either a set of four planes in $\bbP^3$ with an empty intersection or their union. 
We say that each plane entering a tetrahedron $T$ is a face of $T$. 
The intersection line of two faces is its edge, and the intersection point of three faces is its vertex. 

Let us first construct a desmic system of tetrahedra. Fix one tetrahedron $T$.   
Let $\Pi$ be a face of $T$. It defines a projective involution $\sigma_\Pi$ (a harmonic homology) with the set of fixed points equal 
to the union of $\Pi$ and the opposite vertex $v_\Pi$.
Two opposite edges $\ell$ and $\ell'$ (i.e., skew edges) define another 
projective involution $\sigma_\ell$ whose set of fixed points 
is the union of the two edges. It assigns to a point $x$ not on the edges the unique point $x'$ such that 
the pair $\{x,x'\}$ is harmonically conjugate to the intersection points  $\{a,b\}$ 
of the edges with the unique line passing through $x$ and intersecting the edges.

Starting from a general point $P$ not on $T$, we apply the involutions defined by three pairs of opposite edges of $T$. 
Together with $P$, we obtain four points, which are the vertices of a new tetrahedron $T'$. Next, we apply to $P$ four harmonic 
homologies defined by faces of $T$ with its vertices as the centers.  We obtain 
four more points that define a third tetrahedron $T''$. The tetrahedra $T, T'$, and $T''$ obtained in this way are desmic.
Indeed, applying a projective transformation, we may assume that $T = V(xyzw)$ and $P = [1,1,1,1]$. 
The three involutions defined by pairs of opposite edges are given by the scalar matrices with two 
diagonal entries equal to $1$ and two diagonal entries equal to $-1$. So, the vertices of $T'$ are 
$$[1,1,1,1],\ [1,-1,-1,1],\ [1,-1,1,-1],\ [1,1,-1,-1]$$ and 
$$T' = V((x-y-z+w)(x-y+z-w)(x+y-z-w)(x+y+z+w)).$$
 Applying the harmonic homologies to $P$ given by scalar matrices with 
 one diagonal entry equal to $1$ and three diagonal entries equal to $-1$, we find 
that the vertices of $T''$ are 
$$[-1,1,1,1],\ [1,-1,1,1],\ [1,1,-1,1],\ [1,1,1,-1]$$ and 
$$T'' = V((-x+y+z+w)(x-y+z+w)(x+y-z+w)(x+y+z-w)).$$
Next, we check that 
$$-16xyzw+(x-y-z+w)(x-y+z-w)(x+y-z-w)(x+y+z+w)$$
$$+(-x+y+z+w)(x-y+z+w)(x+y-z+w)(x+y+z-w) = 0.$$
This shows that $T, T',$ and $T''$ belong to the same pencil.  

The faces of $T'$ and $T''$ intersect each face of $T$ along the same four lines. 
This gives us 16 lines $V(x,y+z+w), V(x,-y+z+w), V(x,y-z+w),V(x,y+z-w)$, and so on.

It is immediately checked from the coordinates that the three tetrahedra are self-dual with respect to the 
quadric $V(x^2+y^2+z^2+w^2)$. 

\begin{proposition} Three tetrahedra are desmic if and only if each pair is perspective from any vertex of the third one.
\end{proposition}

\begin{proof} Recall that two tetrahedra are perspective from a plane $\Pi$ (resp. from a point $P$) 
if there is a bijection between the sets of their faces (resp. vertices) 
such that $\Pi$ contains four intersection lines of corresponding pairs of faces (resp. the lines joining 
the corresponding vertices intersect at $P$).

Two of the three desmic tetrahedra span a pencil with the base set equal to the union of 16 lines.
They are the intersection lines of $16$ pairs of faces. This implies that each face of one of the tetrahedra contains
four intersection lines of the other two tetrahedra. By definition, each pair of tetrahedra is perspective from each 
face of the third tetrahedron.  By Desargues' Theorem, the two tetrahedra are perspective from each vertex of the 
third tetrahedron (this follows from the proof of 
Desargues's Theorem \cite[Theorem 2.3.1]{CAG}). In fact, since all triples of desmic tetrahedra are projectively 
equivalent (see Proposition \ref{projequiv}), we can check both perspectivity properties from the coordinates. 
Or, check one property and use the self-duality 
of the tetrahedra.

Conversely, the perspectivity of two tetrahedra from each face of the third tetrahedra implies that 
the third tetrahedron contains the base locus of the pencil spanned by the two tetrahedra. Thus, the three 
tetrahedra are desmic.
\end{proof}

\begin{proposition}\label{projequiv} Two desmic systems of tetrahedra are projectively equivalent.
\end{proposition}

\begin{proof} Since two sets of five planes in $\bbP^3$ in general linear position are projectively equivalent, we may assume 
that  $T_1 = V(xyzw)$ is desmic to  $T_2 = V((x+y+z+w)L_1L_2L_3)$.  Since $T_3$ is uniquely determined by $T_1$ 
and $T_2$, it suffices to show that $T_2$ coincides with 
$T'$ from above. The three edges $V(x+y+z+w,L_i)$ of $T_2$ intersect the three pairs of opposite edges of $T_1$.
This gives $L_1 =  a_1(x+y)+b_1(z+w), L_2 =  a_2(x+z)+b_2(y+w), L_3 = a_3(x+w)+b_3(y+z)$. 
Since any edge of $T_2$ intersects a pair of opposite edges of $T_1$,  any pair of the planes $V(L_i)$ 
intersect two edges at the same points. 
For example, $V(L_1)$ intersects the edges of $T_1$ at the points
$$[0,0,1,-1],\  [0,-b_1,0,a_1],\  [0,-b_1,a_1,0],\  [-b_1,0,0,a_1],\   [-b_1,0,a_1,0],\  [-1,1,0,0],$$ 
and $V(L_2)$ intersects the edges of $T_1$ at the points 
$$[0,0,-b_2,a_2],\  [0,1,0,-1],\  [0,-a_2,b_2,0], \ [-b_2,0, 0, a_2], \ [1,0,-1,0],\  [-b_2,a_2,0,0].$$
The common points are $[0,-b_1,a_1,0] =  [0,-a_2,b_2,0]$ and $[-b_1,0,0,a_1] = [-b_2,0, 0, a_2]$. 
This gives $a_1^2 = b_1^2, a_2^2 = b_2^2$. Replacing $L_2$ with $L_3$ we also get $a_3^2=b_3^2$. 
This shows that $T_2$ coincides with $T'$.
\end{proof}

 For example, the desmic pencil
 \begin{equation}\label{eqn1}
 a(x^2+y^2)(z^2+w^2)+b(x^2+w^2)(y^2+z^2)+c(x^2+z^2)(w^2+y^2) = 0,
 \end{equation}
 where $a+b+c = 0$, is projectively equivalent to the desmic pencil spanned by $T,T'$.
 
 Note that the 12 vertices of the three tetrahedra from \eqref{eqn1} lie by pairs on the six edges of each 
 tetrahedron $T, T', T''$, and vice versa.  We say that the two pencils of the desmic tetrahedra are 
 \emph{conjugate } or \emph{associated} (see \cite[Section 3.19]{Lord}).

 \begin{remark}\label{yuz} A pencil of hypersurfaces of degree $d\ge 2$ in $\bbP^n$ 
 contains  
 $k\ge 3$ members equal to the reduced union of $d$ hyperplanes with the empty intersection 
 of any $n+1$ of them 
is called a \emph{desmic pencil}. The desmic pencils in $\bbP^2$ exist for 
 $d\le 5$ \cite{Stipins}. Among them are 
 pencils of conics, the Hesse pencil of cubics, and a desmic pencil 
 of quartics obtained 
 by intersecting a desmic pencil of quartic surfaces with a general plane. A general member of the latter pencil 
 is a smooth plane quartic curve with the property that it admits three bitangents that form the diagonals of a complete 
 quadrilateral, whose six vertices are the points of contact of the bitangents \cite{Humbert2}.

 It is not known whether 
 desmic pencils exist in $\bbP^2$ for $d >5$ and whether there exists a desmic pencil in 
 $\bbP^n, n\ge 3,$ except the desmic pencil of quartics in $\bbP^3$ \cite{Pereira}.
 \end{remark}

\section{Desmic quartic surfaces}
Recall that a \emph{desmic quartic surface} is an irreducible member of the desmic pencil of quartics. 
In fact, any member of the pencil different from the tetrahedra is irreducible. To see this,
 we use the fact that the base locus of the pencil is the union of 16 lines,
no pair of which is skew, lying by four in the faces of any of the tetrahedra. Any reducible member 
must contain an irreducible component of degree one or two. This easily gives a contradiction.

Recall that an \emph{abstract configuration} of type $(a_c,b_d)$ consists of two sets $A$ and $B$ of cardinalities $a,b$ and a relation 
$R\subset A\times B$ such that the fibers of the projection $R\to A$ (resp. $R\to B$) 
have the same cardinality equal to $c$ (resp. $d$). There is a natural definition of an 
isomorphism of abstract configurations. 

In geometry, the relation $R$ is usually the incidence relation between 
points and lines, lines and planes, etc.

\begin{figure}[ht]
\begin{center}
\scalebox{0.3}{
\includegraphics{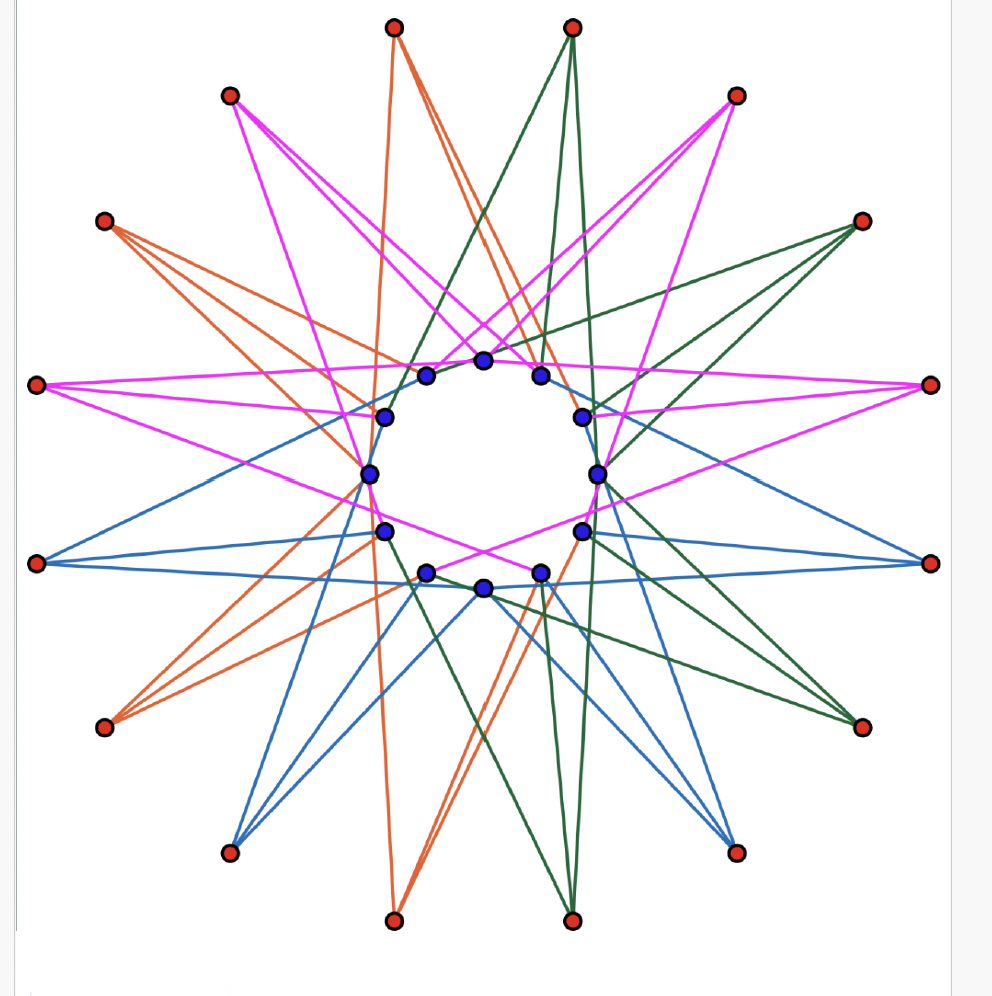}}
\caption{The Levi graph of the Reye configuration (from Wikipedia)}\label{reye}
\end{center}
\end{figure}

The set of 12 singular points and l6 lines on a desmic quartic surface is an example of an abstract configuration 
$(12_4,16_3)$. We call it the \emph{desmic configuration}. Another example of an abstract configuration 
of the same type is the \emph{Reye configuration}. Its set of points consists of eight vertices of a cube in $\bbR^3$, 
the center of the cube, and three 
intersection points at infinity of the 12 edges. Its set of lines consists of 12 edges and 4 diagonals (see Figure \ref{reye}).

\begin{proposition} Let $Q$ be a desmic quartic surface. Then
$Q$ contains $12$ ordinary nodes and $16$ lines that form an abstract configuration $(12_4,16_3)$ isomorphic to the Reye 
configuration. 
\end{proposition} 

\begin{proof} The sixteen lines form the base locus of the pencil. For the pencil from \eqref{eqn1} they are 
$$V(x\pm y,x\pm w),\quad  V(x\pm y,y\pm z),\quad  V(z\pm w,x\pm w),\quad  V(z\pm w, y\pm z).$$
The twelve points are 
\begin{equation}\label{singularpts}
\begin{split}
&[0,0,0,1],\quad [0,0,1,0], \quad [0,1,0,0], \quad [1,0,0,0],\\
&[1,-1,-1,-1],\quad [1,-1,1,1],\quad [1, 1,-1,1],\quad [1,1,1,-1],\\
&[1,1,1,1],\quad [1,1,-1,-1], \quad [1,-1,1,-1],\quad [1,-1,-1,1].
\end{split}
\end{equation}
In the real projective space $\bbP^3(\bbR)$ with the projective coordinates 
$[t_0,t_1,t_2,t_3]$ and the plane at infinity $t_0 = 0$,
the twelve points in the second and the third rows are the vertices of a cube $\Sigma$. The last point in the first row is the center of the cube. 
Each point in the second row is opposite the point directly below it in the third row. The length of an edge of the cube is 
equal to $\sqrt{2}$. The set of 12 points is the set of common vertices of  
the tetrahedra $T, T'$ from above, which are inscribed in the cube. 
The points from the first row form the third inscribed tetrahedron $T''$ in the \emph{complete cube}, that is, the closure $\overline{\Sigma}$ 
of $\Sigma$ in $\bbP^3(\bbR)$.

The singular points lie in pairs at the edges of any tetrahedron. Each face of one of the desmic tetrahedra 
contains four lines.
They form a complete quadrilateral whose three diagonals are the edges lying in the face. 
The incidence property is verified 
immediately. See Figure \ref{quadrilateral} in which six black circles are nodes of $Q$
sitting on the face. 

\begin{figure}[htbp]
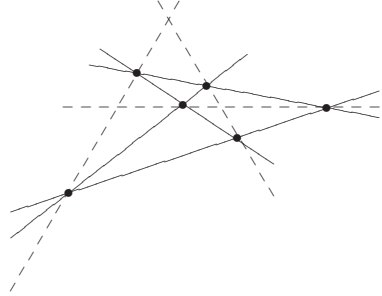

\begin{center}
\scalebox{0.7}{\xy
(-30,0)*{};
(10,0)*{};(70,0)*{}**\dir{--};
(0,-35)*{};(32,20)*{}**\dir{--};
(28,20)*{};(50,-17)*{}**\dir{--};
(0,-25)*{};(45,10)*{}**\dir{-};
(15,8)*{};(70,-2)*{}**\dir{-};
(17,11)*{};(50,-11)*{}**\dir{-};
(70,3)*{};(0,-20)*{}**\dir{-};
@={(43,-6.0),(37.2,3.8),(60,-0.3),(24,6.3),(11,-16.5),(32.7,0.2)}@@{*{\bullet}};
\endxy}
\end{center}
\caption{A complete quadrilateral and three diagonals}
\label{quadrilateral}
\end{figure}
 
Note that the vertices of the desmic tetrahedra 
from \eqref{eqn1} are 
\beq\label{vertices}
[0,0,1,\pm 1],\  [0,1,0,\pm 1], \ [0,1,\pm 1,0], \ [1,0,0,\pm 1], \ [1,0,\pm 1,0], \ [1,\pm 1,0,0]
\eeq
They are different from the singular points of $Q$ and occur as the nodes of a desmic quartic surface from 
the conjugate desmic pencil.  The first six vertices lie in the plane at infinity of the complete 
tetrahedron $\overline{\Sigma}$. They represent the points at infinity of the 12 diagonals of the faces of the cube 
$\Sigma$. The remaining six points are the centers of the faces. Thus, we see that the conjugate desmic systems define two 
Reye configurations originated from the same complete cube.
\end{proof}

Note that the Reye configuration is also isomorphic to the Hesse-Salmon plane configuration of points and lines 
obtained by the projection of the Reye configuration from a general point in the space 
(see \cite[Section 7.3]{DolgachevAbstract}).

\begin{theorem} The group $G$ of the projective transformation that leaves invariant a desmic system of three tetrahedra is 
isomorphic to the quotient of the Weyl group $W(F_4)$ of the root system of type $F_4$ by its center. The group $G$ leaves invariant 
the conjugate desmic system.
\end{theorem}

\begin{proof} Let $\epsilon_1,\ldots, \epsilon_4$ be the standard basis of $\bbR^4$.
The root system $F_4$ consists of $24$ short roots $\pm \epsilon_i$, 
$\pm {1\over 2}(\epsilon_1 \pm \epsilon_2 \pm \epsilon_3 \pm \epsilon_4)$ and
$24$ long roots $\pm (\epsilon_i \pm \epsilon_j)\ (i\not= j)$ \cite[Section 2.10]{Humphreys}. Their images 
in the projective space $\bbP^3(\bbR) = \bbR^4/\{\pm 1\}$ are the vertices of the two conjugate desmic systems of tetrahedra.

The Weyl group $W(F_4)$ is the Coxeter group with the Coxeter diagram (see Figure \ref{coxeter}). 

\begin{figure}[htbp]
\begin{center}
\scalebox{0.9}{\xy
@={(0,0),(10,0),(20,0),(30,0)}@@{*{\bullet}};
(0,0)*{};(30,0)*{}**\dir{-};
(15,3)*{4};
\endxy}
\end{center}
\caption{Coxeter diagram of $W(F_4)$}
\label{coxeter}
\end{figure}
Coxeter's notation $(3,4,3)$ for this group corresponds  to the Schl\"afli symbol $\{3,4,3\}$ of the \emph{$24$-cell}, 
a convex regular polytope in 
$\bbR^4$ with 
$24$ regular octahedra as its boundary (see \cite{Coxeter}). 
The set of long roots can be taken as its set of 
$24$-vertices. 
The dual polytope is also a $24$-cell with the short roots as its vertices. It is known that 
$W(F_4)$ is the group of orthogonal symmetries of the $24$-cell and its dual. The center of $W(F_4)$ generated 
by the element 
$w_0$ of maximal length acts as the minus identity in $\bbR^4$. Passing to the projective space, we see that 
the group $W(F_4)/(w_0)$ acts by projective symmetries of the two conjugate desmic systems. 
Conversely, any projective symmetry 
of the two conjugate desmic systems lifts to an orthogonal symmetry of the $24$-cell. This proves the assertion.
\end{proof}

\vskip6pt

It is known that $W(F_4)$ is isomorphic to the semi-direct product $W(D_4)\rtimes \frakS_3$, where 
$W(D_4)$ is the Weyl group of the root system $R$ of type $D_4$ isomorphic to the semi-direct product 
$2^3\rtimes \frakS_4$ (see \cite[Section 2.12]{Humphreys}, where one can also find 
the relationship with the $24$-cell). The root subsystem of type $D_4$ is generated by the long roots. 
The group $W(F_4)$ leaves this set invariant and acts by orthogonal transformations of the root lattice 
$D_4$ of $R$ generated by the long roots. This defines an injective  homomorphism 
$W(F_4)\to \mathrm{O}(D_4) = W(D_4)\rtimes \frakS_3$, where $\frakS_3$ is the symmetry group of the Dynkin diagram of $D_4$. 
Since both groups have the same order $1152$, $W(F_4) \cong {\mathrm O}(D_4)$.

The desmic pencil from \eqref{eqn1} exhibits an obvious symmetry group isomorphic 
to $2^3\rtimes \frakS_4$ generated by permutations of the coordinates 
$(x,y,z,w)$ and multiplication of coordinates by $\pm 1$. 
The following corollary shows that, in fact, the symmetry group is larger.

\begin{corollary}\label{propnew}
The group $G$ of projective transformations leaving 
a desmic pencil invariant is isomorphic to $W(F_4)/(w_0)$. 
The group of projective symmetries 
acting identically on the base of the pencil is a subgroup of $W(F_4)$ isomorphic to 
$W(D_4)/\{\pm 1\} \cong 2^2\rtimes \frakS_4$. It coincides with the group of projective isomorphisms of any irreducible 
member of the pencil.
\end{corollary}

 \begin{proof} 
 Since all desmic pencils are projectively equivalent, we may assume that our pencil contains the tetrahedra $T, T'$, and $T''$.  
 The projective symmetries of the desmic system obviously leave the desmic pencil invariant. Since any projective symmetry of the pencil 
 leaves the desmic system invariant, we obtain the first assertion of the corollary.
 
 Since the Weyl group acts transitively on the set of its simple roots, it acts transitively on the set of three desmic tetrahedra.
It follows that the image of the group $G$ in its action on the base of the pencil is isomorphic to $\frakS_3$ and 
the kernel of the action is the kernel of the surjective homomorphism $W(F_4)/(w_0)\to \frakS_3$. The element $w_0$ comes 
from the subgroup $W(D_4)$. It is the unique central involution, and it acts as the minus identity on $\bbR^4$. 
A non-trivial automorphism of $\frakS_3$ is bijective. This easily imples that,
under the homomorphism $W(D_4)/(w_0)\rtimes \frakS_3\to \frakS_3$, a subgroup isomorphic to $\frakS_3$ maps isomorphically to 
$\frakS_3$. Thus, the kernel of our action is isomorphic to $W(D_4)/(w_0)$, proving the second assertion.

To prove the last assertion, we use that the projective automorphism group of an irreducible member of the pencil 
leaves invariant the sets of its 12 nodes and 16 lines. It leaves invariant the complete cube $\overline{\Sigma}$ and hence 
leaves invariant the set of three desmic tetrahedra and hence the desmic pencil. 
In the action on the base of the pencil, the group has a fixed point, but the action of $\frakS_3$ on $\bbP^1$ has no 
fixed points for the whole group.
Thus, the projective automorphism group is contained in the group of projective automorphisms of the desmic pencil acting trivially on the its base.
\end{proof}

Note that, in notation of \cite[Table 4.3]{ConwaySmith}, 
the group $W(F_4)$ is isomorphic to 
$\pm \half (\frakS_4\times \frakS_4)\rtimes 2$ and the group $W(F_4)/(w_0) \cong \half (\frakS_4\times \frakS_4)\rtimes 2$.
The latter means that the group is a double extension of the subgroup of $\frakS_4\times \frakS_4$ consisting of pairs
of permutations with the same sign. Although the order $576$ of this group coincides with 
the order of $\frakS_4\times \frakS_4$, the groups are not isomorphic. 
 
 The group $W(F_4)$ also contains a subgroup of index $3$ isomorphic to $W(B_4)\cong 2^4\rtimes \frakS_4$. Its Coxeter's notation is $(3,3,4)$ and it corresponds to 
 the Schl\"afli symbol $\{3,3,4\}$ of the \emph{$16$-cell}, that is, a convex regular polytope in $\bbR^4$ with the boundary equal 
 to the union of $16$ tetrahedra. Its image in $W(F_4)/(w_0)$ is the subgroup $\frac{1}{6}(\frakS_4\times \frakS_4)\rtimes 2$ (in the notation 
 of \cite{ConwaySmith}). This subgroup is realized as the stabilizer subgroup of one of the desmic tetrahedra.

\begin{remark}\label{newrmk} Apparently, Coxeter was the first to relate the Coxeter group $W(F_4)$ with 
the desmic system of tetrahedra (see \cite[Remark 12]{Manivel}). However, the 
relationship between the Reye configuration and the $24$-cell was discussed earlier in 
\cite[Chapter 3, \S\S 22, 23]{Hilbert}.\footnote{The Russian edition has been in the first author's possesion since 1962.} We followed the latter to realize the desmic system inside the 
complete cube. 
\end{remark}

\begin{remark}\label{Szemberg} As observed by Tomasz Szemberg, the set of $24$ points given in \eqref{singularpts} and 
\eqref{vertices} has a 
peculiar property that, although the set is not a complete intersection in $\bbP^3$, its projection to 
$\bbP^2$ from a general point is a complete intersection of a quartic and a sextic curve. The curves are the unions of 
lines. The projection of the six edges of one of the tetrahedra and the projection of the four lines defined by the perspectivity from a vertex of one of the tetrahedra.

\end{remark} 
\vskip6pt

For an elliptic curve $E$, we denote by $\Kum(E\times E)$, called the Kummer surface,
the quotient surface by the inversion of $E\times E$, and by $\widetilde{\Kum}(E\times E)$
the minimal resolution of singularities of $\Kum(E\times E)$.  

\begin{theorem}\label{mainthm1} A desmic surface is birationally isomorphic to the Kummer surface $\Kum(E\times E)$ 
of the self-product of an elliptic curve 
$E$. 
\end{theorem}

\begin{proof}
 Choose two singular points $p_1,p_2$ on two different edges in a face of one of the tetrahedra. Let $E_1,E_2$ be 
 the exceptional curves 
 over these points in the minimal resolution $X$ of the desmic quartic surface $Q$. 
 Among the eight lines that pass through $p_1$ or $p_2$, two coincide with the line $F_0 = \la p_1,p_2\ra$. 
 Let $F_1,\ldots,F_6$ be the remaining 
 lines. Consider the divisor 
 \beq\label{lnsystem}
 H = 2(E_1+E_2)+F_1+\cdots+F_6+2F_0. 
 \eeq
 Since $H\cdot E_i = 1$ and  $H^2 = 4$, the linear system $|H|$ defines a degree
  two map 
  $$\Phi:X \to \bfF_0\subset \bbP^3.$$ 
  Its branch divisor is the union of $8$ lines, 
  four from each ruling. On the other hand, the linear system 
  $|2(E\times \{\textrm{pt}\}+\{\textrm{pt}\}\times E)|$ defines a double cover  
  $\widetilde{\Kum}(E\times E)\to \bbP^1\times \bbP^1$  branched along the union of $8$ lines, from each ruling. 
  Since $\Pic(\bfF_0)$ has no torsion, there is only one isomorphism class of such a cover. This proves 
  the assertion. 
 \end{proof}

\begin{remark}\label{from Kummer}
Conversely, a desmic surface is obtained from the self-product of an elliptic curve $E$ as follows.
Let $\iota:E\to E$ be the inversion map of $E$ and let $a_0=0, a_1, a_2, a_3=a_1+a_2$ be the $2$-torsion points 
on $E$.
Consider the curves $E \times \{a_i\}$, $\{a_i\}\times E$,  
the diagonal $\Delta$, and its translations $\Delta_i$ by a $2$-torsion point 
$\{0\} \times \{a_i\}$ $(i=1,2,3)$, respectively.
A configuration of these curves is given in Figure \ref{16torsions}.

\begin{figure}[htbp]
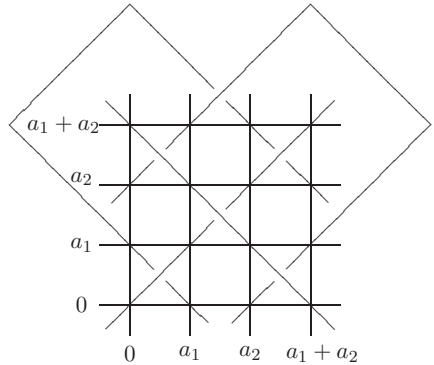

\begin{center}
\scalebox{0.8}{
\xy
(-40,0)*{};
(10,0)*{};(10,40)*{}**\dir{-};
(20,0)*{};(20,40)*{}**\dir{-};
(30,0)*{};(30,40)*{}**\dir{-};
(40,0)*{};(40,40)*{}**\dir{-};
(5,5)*{};(45,5)*{}**\dir{-};
(5,15)*{};(45,15)*{}**\dir{-};
(5,25)*{};(45,25)*{}**\dir{-};
(5,35)*{};(45,35)*{}**\dir{-};
(6,39)*{};(44,1)*{}**\dir{-};
(6,1)*{};(24,19)*{}**\dir{-};
(26,21)*{};(44,39)*{}**\dir{-};
(7,22)*{};(14,29)*{}**\dir{-};
(16,31)*{};(40,55)*{}**\dir{-};
(27,2)*{};(34,9)*{}**\dir{-};
(36,11)*{};(60,35)*{}**\dir{-};
(40,55)*{};(60,35)*{}**\dir{-};
(23,2)*{};(16,9)*{}**\dir{-};
(14,11)*{};(-10,35)*{}**\dir{-};
(-10,35)*{};(10,55)*{}**\dir{-};
(10,55)*{};(24,41)*{}**\dir{-};
(26,39)*{};(34,31)*{}**\dir{-};
(36,29)*{};(43,22)*{}**\dir{-};
(10,-3)*{0};(20,-3)*{a_{1}};(30,-3)*{a_2};(42,-3)*{a_1+a_2};
(2,5)*{0};(2,15)*{a_1};(2,26)*{a_2};(-1,35)*{a_1+a_2};
\endxy
}
\end{center}
\caption{Sixteen 2-torsion points and twelve elliptic curves on $E\times E$}
\label{16torsions}
\end{figure}

The quotient surface of $E\times E$ by the group  $(\bbZ/2\bbZ)^2 = \la \iota \times 1_E, 1_E\times \iota\ra$ is
the quadric $\bbP^1\times \bbP^1$ on which the images of the eight elliptic curves, which are different from the 
four curves
$\Delta$, $\Delta_i$ are 8 lines.  The minimal resolution of the quotient of $E\times E$ by the inversion 
involution $\iota\times \iota$
is $\widetilde{\Kum}(E\times E)$ on which the images of the sixteen $2$-torsion points are sixteen disjoint 
$(-2)$-curves $E_{ij}$ $(0\leq i, j\leq 3)$ and the images of the twelve elliptic curves 
$E\times \{a_i\}$, $\{a_i\}\times E$, $\Delta_i$ are 
twelve disjoint $(-2)$-curves $E_i$, $E^i$, $D_i$ $(0\leq i\leq 3)$.
Let $p_1, p_2$ be the first and the second projections 
$E\times E\to E$, and let
$p_3$ be the map $E\times E \to E$, $(x,y)\to (x+y, x+y)$. 

Then, each $p_i$ induces an elliptic
fibration $f_i: \widetilde{\Kum}(E\times E)\to \bbP^1$ which has four singular fibers of type $\tilde{D}_4$ 
(type ${\rm I}^*_0$ in Kodaira's notation) and eight sections.  
Let $F_i$ be a general fiber of $f_i$ and
let 
\begin{equation}\label{pola}
H = F_1+F_2+F_3 - \half \sum_{0\leq i, j\leq 3}E_{ij}.\footnote{The sum of sixteen disjoint $(-2)$-curves
on a complex $K3$ surface $X$ is divisible by $2$ in $\Pic(X)$ (\cite{Nikulin}).  
In characteristic not equal two, this also holds (see \cite[Remark 3.3]{SB}).}
\end{equation}
Then,
$H^2=4, H\cdot E_{ij}=1$ and $H$ is perpendicular to the remaining twelve $(-2)$-curves.
Thus, the complete linear system $|H|$ gives a birational morphism from $\Kum(E\times E)$ to
a desmic surface on which the images of $E_{ij}$ are sixteen lines and those of the twelve $(-2)$-curves
are twelve nodes.
\end{remark}

The Picard number of $\widetilde{\Kum}(E\times E)$ is at least 19.  For example, consider the linear system
\beq\label{added}
|E_{00}+E_{02}+E_{10}+E_{13}+2(E_0+E_{01}+E^1 + E_{11}+E_1)|,
\eeq
which defines an elliptic fibration with a singular fiber of type $\tilde{D}_8$ and 
sections $E^2$, $E^3$.
The following curves are components of singular fibers of this fibration: $E_{20}$, $E_2$, $E_{22}$, $E_{23}$ 
(forming a Dynkin diagram of type $D_4$) and $D_3$, $E_{30}$, $E_3$, $E_{32}$, $E_{33}$ 
(forming a Dynkin diagram of type $D_5$).  Thus, the Picard lattice contains 
$U\oplus D_8\oplus D_5\oplus D_4$ and hence the Picard number is at least 19.  

 If the Picard number is equal to 19, then it is well known that
the transcendental lattice is isomorphic to $U(2)\oplus \la 4\ra$, where
$\la m\ra$ is the lattice of rank 1 generated by a vector of norm $m$
 (see Remark \ref{transcendental}), which implies that the discriminant of 
$\Pic(\widetilde{\Kum}(E\times E))$ is equal to $\pm 2^4$.  
Since the elliptic fibration defined by the linear system \eqref{added} has 
two sections $E^2, E^3$, $\Pic(\widetilde{\Kum}(E\times E))$ is generated by twenty-eight 
$(-2)$-curves that span an overlattice of $U\oplus D_8\oplus D_5\oplus D_4$ by adding $E^3$, that is, it is 
isomorphic to $U\oplus D_8\oplus D_9$. Thus, we have the following result.

\begin{proposition}\label{PicardLattice} 
The proper transforms of $16$ lines and $12$ exceptional curves of the minimal resolution of a desmic quartic surface 
 forming $(12_4, 16_3)$-configuration generate a quadratic lattice isomorphic to 
 $U\oplus D_8\oplus D_9 (\cong U\oplus E_8\oplus D_8\oplus \la -4\ra)$.
\end{proposition}

\begin{remark}\label{Weierstrass}
In \cite{Humbert}, Humbert showed that the Kummer surface $\Kum(E\times E)$ is birationally
isomorphic to a desmic surface 
by using Weierstrass' $\sigma$-functions (see also \cite[Chapter 2]{Jessop}). Other proofs of this result can be found 
in \cite[Proposition 5.2.3 and Lemma 5.7.6]{Hunt}. Yet another different proof can be found in 
\cite{Manivel}, who relates the desmic pencil with Vinberg's $\theta$-group and its linear representations.

Also, it is known how to compute the $j$-invariant of the elliptic curve $E$ in terms of the 
corresponding parameter of the desmic pencil (see \cite[Theorem 3.10]{HNS}).
\end{remark}

\vskip5pt
\begin{remark} It follows from Corollary \eqref{propnew} that the group of projective
automorphisms of a desmic pencil is isomorphic to $(2^2\rtimes \frakS_4)\rtimes \frakS_3$.
The group acts on the base curve of the pencil. 
The kernel of the corresponding homomorphism is isomorphic to $2^2\rtimes \frakS_4$, 
and the image is isomorphic to $\frakS_3$. 
A desmic quartic with a group of projective automorphisms larger than the group 
$2^2\rtimes \frakS_4$ corresponds to a fixed point 
of the action of $\frakS_3$ on the base. 
Thus, an element of order $3$ corresponds to a point $(a,b,c) = (1,\zeta,\zeta^2)$
where $\zeta^3 = 1, \zeta\ne 1$. 
The surface is isomorphic to $\Kum(E_{\zeta}\times E_{\zeta})$, where 
$E_{\zeta}$ is an elliptic curve with complex multiplication by $\bbZ[\zeta]$. 
Elements of order 2 give three desmic surfaces with $(a,b,c) = (1,1,-2)$ and so on. 
They are isomorphic to $\Kum(E_{\sqrt{-1}}\times E_{\sqrt{-1}})$, where 
$E_{\sqrt{-1}}$ is an elliptic curve with complex multiplication by $\bbZ[\sqrt{-1}]$.
For the (full) automorphism groups of these Kummer surfaces, see \cite[Theorems 5.4, 5.6, 5.7]{KeKo}.

The group $\frakS_3$ acts on the desmic pencil permuting the desmic tetrahedra. 
The action has two fixed points with the stabilizer subgroup of order $3$ 
and three fixed points with the stabilizer subgroup 
of order $2$. This defines a Galois $\frakS_3$-cover from the desmic pencil to $\bbP^1$. 
As is well known, it is isomorphic to the modular curve $X(\Gamma_0(2))$, 
and the cover is isomorphic to the natural projection to the modular curve 
$X(\Gamma(1))\cong \bbP^1$. This agrees with \cite[Section B.5, Corollary B.5.7]{Hunt}. 

It follows 
from Remark \ref{from Kummer} that the moduli space of desmic K3 surfaces is isomorphic to 
the moduli space of Kummer surfaces $\widetilde{\Kum}(E\times E)$ together 
with a choice of one of the three projections $E\times E\to E$. 
This shows that the moduli space of desmic K3 surfaces is isomorphic to 
a $\frakS_3$-cover of the moduli space of elliptic curves isomorphic to $X(\Gamma(1))$. 
\end{remark}

\begin{remark}\label{transcendental} 
Let $T(E\times E)$ be the transcendental lattice of the abelian surface $E\times E$. 
 As is well known, this lattice is isomorphic to $U\oplus \la 2\ra$ for a generic $E$. The transcendental 
 lattice of the minimal nonsingular model $\widetilde{\Kum}(E\times E)$ is isomorphic 
 to $T(E\times E)(2)$. Hence, the Picard lattice of $\widetilde{\Kum}(E\times E)$ is isomorphic to the lattice 
 $M_2 = U\oplus E_8\oplus D_8\oplus \la -4\ra$. The special cases of curves $E_\omega$ with complex multiplication 
 have the Picard lattices 
 isomorphic to $U\oplus E_8\oplus E_8\oplus \la -4\ra\oplus \la -4\ra$ if $\omega = \sqrt{-1}$ and 
 $U\oplus E_8\oplus E_8\oplus A_2(2)$ if 
 $\omega = \zeta$.
 
 \end{remark}

\section{The cubic complex of lines}
Let $\calN$ be a net of quadrics in $\bbP^3$. We assume that $\calN$ contains a smooth quadric. The set of lines 
contained in some quadric from $\calN$ form  
a \emph{line complex} $\frakG(\calN)$, a hypersurface
in the Grassmannian $G_1(\bbP^3)$ of lines in $\bbP^3$. It is called the Montesano complex, and its 
degree is equal to three. Recall that the latter means that $\frakG(\calN)$ is a complete intersection of 
$G_1(\bbP^3)$, viewed as a quadric in the Pl\"ucker embedding in $\bbP^5$, with a cubic hypersurface.

The degree of a line complex $\frakG$ is equal to the degree of the plane curve 
$\Omega(x)\cap \frakG$, where $\Omega(x)$ is the plane in $G_1(\bbP^3)$ of lines 
containing a general point $x\in \bbP^3$. We easily see why the degree of $\frakG(\calN)$ is equal to $3$.
The set of all quadrics from $\calN$ that contain $x$ is a pencil $\calN(x)$ whose base curve is a quartic 
elliptic curve $C$ containing $x$.
The lines through $x$ contained in a quadric from $\calN(x)$ form a cone over the projection of $C$ from $x$.
It is a cubic cone, and $\Omega(x)\cap \frakG(\calN)$ is the projection of $C$, a cubic curve.

Let us return to our desmic pencil given in (\ref{eqn1}). Take two of the desmic tetrahedra, say 
$V((x^2-y^2)(z^2-w^2))$ and $V((x^2-z^2)(y^2-w^2))$. The quadrics 
$V((x-y)(z-w)), V((x-z)(y-w)), V((x+w)(y+z))$ vanish at $8$ points given in the first two rows in \eqref{singularpts}.
It is immediately checked 
that the quadrics are linearly independent and span a net $\calN_1$ of quadrics with $8$ base points. 
It obviously contains a smooth quadric. The discriminant curve of $\calN_1$ is a plane quartic curve in $\calN_1$ with 
three nodes corresponding to three reducible quadrics from the net. Similarly, we find the nets $\calN_2, \calN_3$ 
with base points 
given by the first and the third rows, and by the second and the third rows:
$$\calN_2= \la V((x-y)(z+w)), V((x-z)(y+w)), V((x-w)(y+z))\ra;$$
$$\calN_3=\la V(x^2-y^2), V(x^2-z^2), V(x^2-w^2)\ra.$$
Thus, a desmic pencil defines three nets of quadrics $\calN_1, \calN_2, \calN_3$, and hence, 
three cubic line complexes
$\frakG_1, \frakG_2,\frakG_3$.

\begin{corollary} The $12$ lines that join a point $p\in \bbP^3$ with the $12$ nodes 
of a desmic quartic $Q$ are contained in a cubic cone with its vertex at $p$.
\end{corollary}

\begin{proof} If $p$ is a base point of $\calN_1$, the polar cubic surface 
with the pole at $p$ is a cubic cone that contains all singular points of $Q$. So, we may assume that $p$ is not a 
base point.
Then, a line spanned by $p$ and one of 
the base points of $\calN_1$ is contained in a quadric from $\calN_1$. Hence, it belongs to 
the cubic line complex, and, therefore, belongs to the cubic cone $K_1(p)$ of lines from $\frakG_1$ passing 
through the point $p$. This shows that the ruling of the cubic cone $K_1(p)$ contains lines spanned by $p$ and 
any of the eight base points of $\calN_1$. Let $p_i$ be one of the remaining four nodes of $Q$.
There are four lines joining $p_i$ with two base points of $\calN_1$. Thus, the pencil 
$\calN_1(p_i)$ of quadrics from 
$\calN_1$ passing through $p_i$ contains four lines passing through $p_i$. 

The lines through $p_i$ contained in a quadric from $\calN_1(p_i)$
must be a cone with its vertex at $p_i$.
So, we choose the quadric from this pencil passing through $p$ and obtain that the line $\la p,p_i\ra$ belongs to 
the cubic line complex $\frakG_1$. 
Thus, the line $\la p,p_i\ra$ is a generator of the cubic cone $K_1(p)$.
\end{proof}

Since the cubic cones $K_i(x)$ of different line complexes $\frakG_i$ intersect along 12 lines, they coincide.
This gives the following corollary.

\begin{corollary} The cubic complexes $\frakG_1, \frakG_2, \frakG_3$ coincide.
\end{corollary}

Let $\frakG = \frakG_1= \frakG_2 = \frakG_3.$ 

\begin{defn} The cubic line complex  $\frakG$ arising from a desmic pencil is called the \emph{Humbert desmic line complex}.
\end{defn}

Since $\frakG$ contains three plane-pairs formed by faces of the desmic tetrahedra, the corollary implies 
another corollary.

\begin{corollary} Any line contained in the face of one of the desmic tetrahedra belongs to $\frakG$.
\end{corollary}

\begin{remark} Since any of the twelve nodes of the desmic quartic $Q$ is a base point 
of one of the nets $\calN_i$, any line through a node belongs to the complex 
$\frakG$. Thus, $\frakG$ contains twelve $\alpha$-planes $\Lambda_i$ of lines 
passing through a node and twelve $\beta$-planes $\Xi_i$ of lines contained 
in some face of one of the desmic tetrahedra. 
Each of the 16 lines $\ell_i$ in $Q$ is common to three planes $\Lambda_i$ and three planes $\Xi_i$. 
Since the six planes span the whole space $\bbP^5$, 
the 16 lines are singular points of $\frakG$. The 16 singular points and 12 planes from each family of planes in 
the Grassmannian form two configurations $(12_4,16_3)$. Together, they form a configuration 
$(16_6,24_4)$. 

Note that each edge of a desmic tetrahedron is contained in two 
$\alpha$-planes and two $\beta$-planes. Since the four planes span the whole space $\bbP^5$, 
each edge is a singular point of $\frakG$. Therefore, 
the cubic line complex contains $16+18 = 34$ singular points. Each $\beta$-plane contains seven singular points 
of $\frakG$; they form a 
complete quadrilateral and its three diagonals. Each $\alpha$-plane contains five singular points, one edge, 
and two lines $\ell_i$ from each face containing 
the edge.

The group $(2^2\rtimes \frakS_4)\rtimes \frakS_3$ of symmetries of the desmic pencil 
acts on $\frakG$. This makes $\frakG$ a Fano $3$-fold similar 
to the Segre cubic primal that contains 10 singular points and 15 planes forming 
an abstract configuration 
$(15_4,10_6)$ with the group of symmetries isomorphic to $\frakS_6$. 
\end{remark}

\vskip5pt
\begin{lemma} Each quadric from $\calN_i$ cuts a desmic quartic surface $Q$ along the union of two quartic curves 
intersecting at 
the base points of $\calN_i$.
\end{lemma} 

\begin{proof} A quadric  $G\in \calN_i$ intersects $Q$ along a curve $C$ of degree eight with $8$ double points. 
A general quadric from the pencil of quadrics in $\calN_i$ passing through a general point $x\in C$ intersects $C$ with 
multiplicity $2\times 8+1 = 17 > 16$. This shows that the base curve of the pencil $\calN_i(x)$ is a part of $C$.
The residual part is another quartic curve.
\end{proof}

The next proposition provides an alternative proof that a desmic quartic surface is birationally isomorphic to $\Kum(E\times E)$.

\begin{proposition} Let $X$ be a minimal resolution of singularities of $Q$. 
Then each net of quadrics $\calN_i$ defines 
a base-point-free pencil  
 $|E_i| = |E_i'|$ of elliptic curves on $X$, where $E_i,E_i'$ are the proper transforms of the irreducible components 
 of the intersection of a quadric from $\calN_i$ with $Q$.  We have 
 $E_i\cdot E_j = 2$, and the linear system $|E_i+E_j|$ defines a degree two cover  $X\to \bbP^1\times \bbP^1$ with 
 the branch curve 
 equal to the union of 8 lines, four from each family of lines on the quadric model of $\bbP^1\times \bbP^1$.
 \end{proposition}

This easily follows from the discussion above. 

\begin{proposition}\label{16conics} A desmic quartic surface contains $16$ conics. 
\end{proposition}

\begin{proof} Let $\Pi$ be a plane containing one of the sixteen lines $\ell$ on $Q$. For example, we may take 
$\ell$ to be $ V(x+y,z+w)$. Then, $\Pi = V(u(x+y)+v(x+w))$. Substituting in the equation of $Q$, we find 
that $\Pi$ is tangent to $Q$ along $\ell$ provided that $u(b+c)+v(a+b) = 0$. The residual intersection of $Q\cap \Pi$ 
is a 
conic intersecting $\ell$ at two different points. 
\end{proof}

It follows from the proof of Proposition \ref{16conics} that the pencil of cubic curves cut out by planes 
through any of the lines contains three singular fibers of type $\tilde{A}_5$, and one singular fiber of type 
$\tilde{A}_1$. The fibers of type $\tilde{A}_5$ 
come from the three quadrilaterals containg $\ell$, and the fiber of type $\tilde{A}_1$ comes from the union of $\ell$ 
and the residual conic.
 
 \vskip5pt 
Let us find the equation of the Humbert desmic cubic line complex $\frakG$.
The net of quadrics passing through the eight points given in rows one and two in \eqref{singularpts} is equal 
to 
$$t_1(xy+zw)+t_2(xz+yw)+t_3(xw+yz) = 0.$$
Let 
$$(x,y,z,w) = (a_1u+b_1v,a_2u+b_2v,a_3u+b_3v,a_4u+b_4v)$$
be the parametric equation of a line in $\bbP^3$. Plugging it in the previous equation, we obtain that condition 
that there exist parameters $(t_1,t_2,t_3)$ such that the quadric contains the line is 
\[
\det
\scalebox{0.8}
{$\begin{pmatrix}
a_1a_2+a_3a_4&a_1a_3+a_2a_4&a_1a_4+a_2a_3\\
b_1b_2+b_3b_4&b_1b_3+b_2b_4&b_1b_4+b_2b_3\\
a_1b_2 + a_2b_1 + a_3b_4 + a_4b_3&a_1b_3 + a_2b_4 + a_3b_1 + a_4b_2&a_1b_4+a_2b_3+a_3b_2+a_4b_1
\end{pmatrix}$}
=0.
\]
Computing the determinant, we obtain that the equation of $\frakG$ in Pl\"ucker coordinates $(x_1,x_2,x_3,x_4,x_5,x_6)=(p_{12},p_{13},p_{14},p_{23},p_{24},p_{34})$ is the following:
\begin{eqnarray}\label{eqCubic2}
&p_{12}p_{34} - p_{13}p_{24} + p_{14}p_{23} = 0,\\ \label{eqCubic3}
&-p_{12}p_{13}p_{23} + p_{12}p_{14}p_{24} - p_{13}p_{14}p_{34} + p_{23}p_{24}p_{34} = 0. 
\end{eqnarray}
Note that the equation is invariant with respect to the natural action of $\frakS_4$ tensored 
with the sign representation. It is also invariant with respect to the action of $2^3$ on coordinates in $\bbP^3$.
Thus, we obtain that the group $2^3\rtimes \frakS_4$ acts on $\frakG$. This agrees with the symmetry group of the 
desmic pencil.
The set $\Sing(\frakG) = \Sing(\frakG)_1\cup \Sing(\frakG)_2$ of 34 nodes of 
$\frakG$, where the subset $ \Sing(\frakG)_1$ consists of 18 points

\begin{equation}\label{18points}
\def\no#1{{\scriptstyle#1}\colon}
\def\sp{\ \,}
\scalebox{0.8}{$\alignedat6
 \no{1}[1,0,0,0,0,0],\ \no{2}[0,1,0,0,0,0],\  \no{3}[0,0,1,0,0,0],\quad  \qquad \qquad\\
 \no{4}[0,0,0,1,0,0],\ \no5[0,0,0,0,1,0],\ \no6[0,0,0,0,0,1],\quad \qquad \qquad\\
\no7[1,1,0,0,1,1],\no{8}[1,1,0,0,-1,-1],
\no9[1,-1,0,0,1,-1],\no{10}[1,-1,0,0 ,-1,1],\\
 \no{11}[1,0,1,1,0,-1],\no{12}
[1,0,1,-1,0,1],\ 
\no{13}[1,0,-1,1,0,1],\ \no{14} [1,0,-1,-1,0,-1],\\
 \no{15}[0,1,1,1,1,0], 
\no{16}[0,1,1,-1,-1,0], \no{17}[0,1,-1,1,-1,0],\no{18}[0,1,-1,-1,1,0].
\endalignedat$}
\end{equation}
corresponding to 18 edges which form two orbits of cardinalities 6 and 12, and
$\Sing(\frakG)_2$ consists of 16 points
\begin{equation}\label{16points}
\def\no#1{{\scriptstyle#1}\colon}
\def\sp{\ \,}
\scalebox{0.8}{$\alignedat4
\no1&[1,1,1,0,0,0], \sp&\no2& [1,1,-1,0,0,0], \sp&\no3 & [1,-1,1,0,0,0], \sp&\no4 &[1,-1,-1,0,0,0]\\
\no5&[1,0,0,1,1,0], \sp &\no6& [1,0,0,1,-1,0], \sp&\no7 & [1,0,0,-1,1,0],\sp&\no8 & [1,0,0,-1,-1,0]\\
\no9&[0,1,0,1,0,1], \sp &\no{10}& [0,1,0,1,0,-1], \sp&\no{11}  & [0,1,0,-1,0,1], \sp &\no{12}& [0,1,0,-1,0,-1]\\
\no{13}&[0,0,1,0,1,1], \sp&\no{14} & [0,0,1,0,1,-1], \sp &\no{15}& [0,0,1,0,-1,1],\sp&\no{16}  & [0,0,1,0,-1,-1]
\endalignedat$}
\end{equation}
corresponding to 16 lines.

The set of 24 planes contained in $\frakG$ consists of 12 $\alpha$-planes and 12 $\beta$-planes.
Recall that the Pl\"ucker equations of the $\alpha$-plane of lines passing 
through a fixed point $[a,b,c,d]\in \bbP^3$ are 
$$-cp_{12}+bp_{13}-ap_{23} = dp_{13}-cp_{14}+ap_{34} = dp_{12}-bp_{14}+ap_{24} = 0.$$
The equations of the $\beta$-plane of lines in a fixed plane $V(ax+by+cz+dw)\subset \bbP^3$ are 
\beq
\begin{split}
&cp_{12}-bp_{13}+ap_{24} = dp_{13}-bp_{14}+ap_{24} = -dp_{13}+cp_{14}+ap_{34} = 0, \quad a\ne 0,\\ 
&-cp_{12}+bp_{13} = -dp_{12}+bp_{14} = dp_{23}-cp_{24}+bp_{34} = 0, \quad a = 0.
\end{split}
\eeq

For brevity of the notation, let $(p_{12},p_{13},p_{14},p_{23},p_{24},p_{34}) = 
(x_1,x_2,x_3,x_4,x_5,x_6)$. Substituting the coordinates of points from \eqref{singularpts}, 
we obtain the equations of the twelve 
 $\alpha$-planes:
\begin{equation*}
\def\no#1{{\scriptstyle#1}\colon}
\scalebox{0.8}{$\alignedat3
\no1&V(x_1, x_2, x_4),\no{2} V(x_1, x_3, x_5),
\no3 V(x_2, x_3, x_6),\no4 V(x_4, x_5, x_6),\\
\no{5/6}&V(x_1+ x_2 + x_4, \pm x_1 + x_3+x_5, \pm x_2 - x_3+x_6),
\no{7/8}V(x_1+ x_2 - x_4, \pm x_1- x_3+ x_5, \pm x_2 +x_3+ x_6),\\
 \no{9/10}&V(x_1- x_2 + x_4, \pm x_1 -  x_3 + x_5, \pm x_2- x_3+ x_6),
\no{11/12}V(x_1- x_2 - x_4, \pm x_1+ x_3 + x_5, \pm x_2+ x_3+ x_6).
\endalignedat$}
\end{equation*}
Substituting the equations of the faces of the three tetrahedra from \eqref{eqn1}, we obtain the equations 
of the twelve $\beta$-planes:
\begin{equation*}
\def\no#1{{\scriptstyle#1}\colon}
\scalebox{0.8}{$\alignedat2
\no{1/2}V(x_1,x_2\pm x_4, x_3\pm x_5),\no{3/4}V(x_2, x_1\pm x_4, x_3\mp x_6),\no{5/6} 
V(x_3, x_1\pm x_5, x_2\pm x_6),\\
\no{7/8}V(x_4, x_1\pm x_2, x_5\pm x_6),\no{9/10} V(x_5, x_1\pm x_3, x_4\mp x_6),\no{11/12} V(x_6, x_2\pm x_3, x_4\pm x_5).
\endalignedat$}
\end{equation*}

\vskip5pt
 Changing the Pl\"ucker coordinates to the Klein coordinates 
$$(x_1,x_2,x_3,y_1,y_2,y_3)$$
$$ = (p_{12}+p_{34},\ -p_{13}+p_{24},\ p_{14}+p_{23},\
i(p_{34}-p_{12}),\  i(p_{24}+p_{13}),\  i(p_{23}-p_{14})),$$
where $i=\sqrt{-1}$,
we obtain the new equation of $\frakG$:
\beq\label{humbertcomplex}
x_1^2+x_2^2+x_3^2+y_1^2+y_2^2+y_3^2 = x_1x_2x_3+iy_1y_2y_3= 0.
\eeq
It is in striking resemblance to Baker's equations of the Segre cubic primal:
$$x_1+x_2+x_3+y_1+y_2+y_3 = x_1x_2x_3+y_1y_2y_3 = 0$$
(see \cite[Chapter V, Ex. 20]{BakerBook}).
Note that the presence of $\sqrt{-1}$ in the equation is very important. 
Replacing $\sqrt{-1}$ with $1$, we obtain a threefold with only 18 nodes.

We can rewrite the coordinates of the singular points and the equations of the 24 planes in the new coordinate system.

The set $\Sing(\frakG) = \Sing(\frakG)_1\cup \Sing(\frakG)_2$, where 
the subset $ \Sing(\frakG)_1$ consists of 18 points
\begin{equation}\label{18newpoints}
\def\no#1{{\scriptstyle#1}\colon}
\def\sp{\ \,}
\scalebox{0.8}{$\alignedat6
 \no{1}[1,0,0,-i,0,0],\ \no{2}[0,1,0,0,-i,0],\  \no{3}[0,0,1,0,0,-i],\quad  \qquad \qquad\\
 \no{4}[0,0,1,0,0,i],\ \no5[0,1,0,0,i,0],\ \no6[1,0,0,i,0,0],\quad \qquad \qquad\\
\no7[1,0,0,0,i,0],\no{8}[0,1,0,i,0,0],
\no9[0,1,0,-i,0,0],\no{10}[1,0,0,0 ,-i,0],\\
 \no{11}[0,0,1,-i,0,0],\no{12}
[1,0,0,0,0,-i],\ 
\no{13}[1,0,0,0,0,i],\ \no{14} [0,0,1,i,0,0],\\
 \no{15}[0,0,1,0,i,0], 
\no{16}[0,1,0,0,0,i], \no{17}[0,1,0,0,0,-i],\no{18}[0,0,1,0,-i,0].
\endalignedat$}
\end{equation}
and $\Sing(\frakG)_2$ consists of $16$ points
\beq\label{16newpoints}
[\epsilon_1,\epsilon_2, \epsilon_3,e_1,e_2,e_3],
\eeq
where $\epsilon_i^2 = -1, e_i^2 = 1,$ and 
$\epsilon_1\epsilon_2\epsilon_3+ie_1e_2e_3 = 0$.

The equations of planes become 
\beq\label{neweqplanes}
V(x_1-\epsilon_1\sqrt{-1}y_{\sigma(1)},x_2-\epsilon_2\sqrt{-1}y_{\sigma(2)}, x_3-\epsilon_3\sqrt{-1}y_{\sigma(3)}),
\eeq
where $\epsilon_i = \pm 1$, and $\epsilon_1\epsilon_2\epsilon_3 = 1$.

We can identify the equations with permutations  
$g\in \frakS_4$ as follows. The vector 
$(\epsilon_1,\epsilon_2,\epsilon_3)$ is identified 
with elements from the normal subgroup $2^2$ of $\frakS_4$. 
For example, if we 
let $ a = (14)(23), b = (13)(24), c = (12)(34) $,
 then
$$(1,1,1)\leftrightarrow 1,\  (1,-1,- 1)\leftrightarrow a, \ (-1,1,- 1)\leftrightarrow b, 
(-1,-1,1)\leftrightarrow c.$$
Now,  we can 
write elements of $\frakS_4$ as the products of elements from $2^2$ and elements from 
$\frakS_3$ generated by 
$(12), (23)$. The elements from $\frakS_4\setminus \frakS_3$ are
\[
{\small \begin{split}
&(14) = (23)a,\  (24) = (13)b,\  (34) = (12)c,\ 
 (1234) = (24)a, \ (1243) = (23)b,\\
 &(1342) = (23)c,\ 
(1324) = (34)a,\  (1423) = (34)b,\  (1432)= (24)c,\\
&(124) = (132)b,\  (142) = (123)a,\  (234) = (132)c,\\
& (243) = (123)b,\  (134) = (123)c,\  (143) = (132)a.
\end{split}}
\]

The set of $\alpha$-planes now becomes
\begin{equation*}
\def\no#1{{\scriptstyle#1}\colon}
\scalebox{0.9}{$\alignedat6
\no1&(12)(34),\quad \no{2}&(13)(24),\quad 
\no3& (14)(23),\quad \no4& 1,\quad
\no{5}&(142), \quad \no{6}&(132),\\
\no{7}&(123),\quad \quad  \no8&(124),\quad 
 \no{9}&(143),\quad  \no{10}&(243),\quad  \no{11}&(234),\quad  \no{12}&(134).
\endalignedat$}
\end{equation*}

The set of $\beta$-planes becomes:
\begin{equation*}
\def\no#1{{\scriptstyle#1}\colon}
\scalebox{0.9}{$\alignedat6
\no1&(1342),\quad  \no{2}&(1243),\quad  
\no3&(1432),\quad  \no4&(1234),\quad 
 \no{5}&(1423),\quad  \no{6}&(1324),\\
\no{7}&(12),\quad \quad \no8&(34),\quad  
 \no{9}&(24),\quad \quad \no{10}&(13),\quad 
\no{11}&(23),\quad \quad \no{12}&(14).
\endalignedat$}
\end{equation*}

\vskip5pt
Let us examine the incidence relation between $34$ singular points and $24$ planes.

Let $H_1, H_2, H_3$ be subgroups of
$\frakS_4$ generated by two commutating  transpositions 
$\{(12), (34)\}$, $\{(13), (24)\}$, $\{(14), (23)\}$, respectively.  
Then we can identify 18 left
cosets of $H_1$, $H_2$, $H_3$ and 18 points from $\Sing(\frakG)_1$, and   
the incidence graph of $(24_3,18_4)$ is the same as that of 24 elements of $\frakS_4$ 
and 18 cosets. 

 For example, the point $(i,0,0,0,0,1)$ from $\Sing(\frakG)_1$ is contained in four planes
$V(x_1-iy_3,x_2\pm iy_1,x_3\pm iy_2)$, and $V(x_1-iy_3,x_2\pm iy_2,x_3\pm iy_1)$. They correspond to the permutations
$(143),(132)$ and $(1432),(13)$, respectively. The four permutations form one coset 
of the subgroup $H_1$.

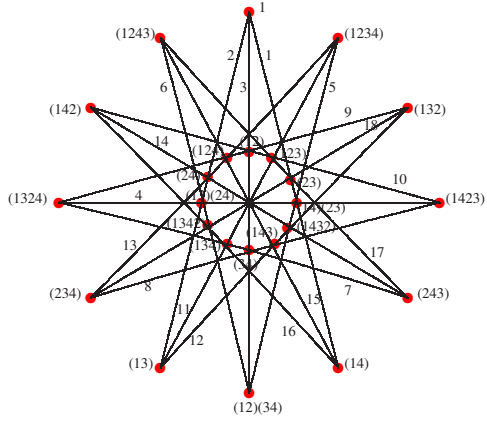
\begin{figure}[ht]
\begin{center}
\unitlength1.5pt
\def\r{\color{blue}\circle*{3}}
\def\b{\color{red}\circle*{3}}
\def\c{\color{blue}\circle*{3}}
\scalebox{0.8}{\begin{picture}(60,60)(-20,10)
\put(0,60)\b  \put(3,60){\tiny{1}}
\put(0,-60)\b \put(-5,-66){\tiny{(12)(34)}}
\put(15,0)\b \put(15,-3){\tiny{(14)(23)}}
\put(-15,0)\b \put(-20,1){\tiny{(13)(24)}}
\put(-28,-52)\b  \put (-38,-52){\tiny{(13)}}
\put(60,0)\b  \put(62,0){\tiny{(1423)}}
\put(0,-15)\b  \put(-5,-21){\tiny{(34)}} 
  \put(-13,8)\b \put(-23,7){\tiny{(24)}}
\put(28,-52)\b \put (30,-52){\tiny{(14)}}
\put(-13,-7)\b  \put(-26,-8){\tiny{(1342)}}
\put(7,14)\b  \put(8,14){\tiny{(123)}}
\put(-50,30)\b  \put(-63,28){\tiny{(142)}}
\put(50,30)\b  \put (52,28){\tiny{(132)}}
\put(-50,-30)\b  \put (-63,-30){\tiny{(234)}}
\put(50,-30)\b  \put (53,-30){\tiny{(243)}}
\put(8,-13)\b  \put(-1,-11){\tiny{(143)}}
\put(-7,-13)\b \put(-19,-15){\tiny{(134)}}
\put(-7,14)\b  \put(-18,15){\tiny{(124)}}
\put(0,16)\b \put(-3,18){\tiny{(12)}}
\put(13,7)\b \put(15,5){\tiny{(23)}}
\put(-60,0)\b  \put(-76,0){\tiny{(1324)}}
\put(28,52)\b  \put (30,52){\tiny{(1234)}}
\put(12,-8)\b  \put(15,-9){\tiny{(1432)}} 
\put(-28,52)\b  \put (-42,52){\tiny{(1243)}}
\put(-3,35){\tiny{3}}
\put(-7,45){\tiny{2}}
\put(5,45){\tiny{1}}
\put(-28,35){\tiny{6}}
\put(25,35){\tiny{5}}
\put(-36,1){\tiny{4}}
\put(18,-32){\tiny{15}}
\put(10,-42){\tiny{16}}
\put(36,23){\tiny{18}}
\put(38,-17){\tiny{17}}
\put(-23,-35){\tiny{11}}
\put(-19,-45){\tiny{12}}
\put(-40,-15){\tiny{13}}
\put(-30,18){\tiny{14}}
\put(30,27){\tiny{9}}
\put(45,6){\tiny{10}}
\put(30,-29){\tiny{7}}
\put(-33,-28){\tiny{8}}
\qbezier(0,60)(0,0)(0,-60) 
\qbezier(0,60)(-14,4)(-28,-52) 
\qbezier(0,60)(14,4)(28,-52) 
\qbezier(28,-52)(0,0)(-28,52)
\qbezier(0,-60)(-14,0)(-28,52)
\qbezier(50,-30)(10,10)(-28,52)
\qbezier(28,52)(14,0)(0,-60) 
\qbezier(28,52)(-10,10)(-50,-30)
\qbezier(28,52)(0,0)(-28,-52)
\qbezier(50,30)(12,-8)(-28,-52)
\qbezier(50,30)(0,0)(-50,-30)
\qbezier(50,30)(-5,15)(-60,0)\
\qbezier(60,0)(5,15)(-50,30)
\qbezier(60,0)(0,0)(-60,0)
\qbezier(60,0)(-1,-15)(-50,-30)
\qbezier(50,-30)(0,0)(-50,30)
\qbezier(50,-30)(5,-15)(-60,0)
\qbezier(28,-52)(-10,-10)(-50,30)
\end{picture}}
\end{center}
\vskip80pt
\caption{The configuration $(24_3,18_4)$}
\label{FigureC3}
\end{figure}

The abstract configuration $(24_3,18_4)$ defined by the incidence relation between planes and points is the 
configuration $C_3$ from \cite[Figure 1.7.4]{Grunbaum} as in Figure \ref{FigureC3}.

\vskip20pt

\vskip20pt
The incidence graph of the configuration $(24_4,16_6)$ coincides with 
the incidence graph of elements of a $4\times 4$ matrix $(a_{ij})$ and 
monomials $a_{1j_1}a_{2j_2}a_{3j_3}a_{4j_4}$ entering 
the expression of the determinant of the matrix.  For example, let us identify $(a_{ij})$
with $\Sing(\frakG)_2$ given in (\ref{16points}) as follows:
\begin{equation*}
\def\no#1{{\scriptstyle#1}}
(a_{ij}) = 
\scalebox{0.8}{$\begin{pmatrix}
\no1&\no{14} & \no{12} & \no7\\
\no{15} & \no2 &\no5 &\no{10} \\
\no9& \no8&\no3 & \no{16} \\
\no6&\no{11}& \no{13}& \no4 \\
\end{pmatrix}$}.
\end{equation*}
It is immediately checked that  
each plane contains $4$ points from $\Sing(\frakG)_2$ and each point from $\Sing(\frakG)_2$ is contained in six planes.
We say that the configuration of planes and points is of type 
$(24_{3+4}, 18_4+16_6)$ meaning that each plane contains three points from the set 
$\Sing(\frakG)_1$ 
and $4$ points from the set $\Sing(\frakG)_2$. Also, each point from $\Sing(\frakG)_1$ is contained in four planes, 
and each point 
of the set $\Sing(\frakG)_2$ is contained in $6$ planes.

The configuration $(24_{3+4}, 18_4+16_6)$ is the union of two configurations 
$(24_3,18_4)$ and $(24_4,16_6)$. In both configurations, we can identify 
the set of 24 elements with the set of permutations $\frakS_4$.

\vskip6pt

 The new coordinates with the equation \eqref{humbertcomplex} exhibit  an obvious group of projective symmetries of $\frakG$ isomorphic to 
 $2^4\rtimes (\frakS_3\times \frakS_3)\cong \frakS_4\times \frakS_4$.
 
 \begin{proposition}\label{symmetries} Assume $p\ne 2$. The group of projective automorphisms of the cubic line 
 complex $\frakG$ is isomorphic 
to the group $G: = \frakS_4\wr 2$. The additional generator is a projective automorphism of 
order $2$ defined by 
$$g_0:(x_1,x_2,x_3,y_1,y_2,y_3)\mapsto (-y_1,-y_2,-y_3,x_1,x_2,x_3).$$ 
 The group $G$ has two orbits 
$\Sing(\frakG)_1$ and $\Sing(\frakG)_2$ on the set $\Sing(\frakG)$ and one orbit on the set $\calP(\frakG)$ of planes 
in $\frakG$.
\end{proposition}
 \begin{proof} We immediately check that $g_0$ is a projective automorphism of $\frakG$ of order $2$ that normalizes 
 the subgroup $\frakS_4\times \frakS_4$, and that the subgroup 
$G = (\frakS_4\times \frakS_4)\rtimes \la g_0\ra$ of $\Aut_{\textrm{proj}}(\frakG)$ 
has two orbits on the set $\Sing(\frakG)$ and one orbit on the set $\calP(\frakG)$ of planes.

Let $G_1 =  \Aut_{\textrm{proj}}(\frakG)$. 
Since $G$ acts transitively and faithfully on the 
set $\calP(\frakG)$, it suffices to see that the stabilizer subgroup $H\subset G_1$ of a plane $\Pi\in \calP(\frakG)$ 
is isomorphic to $\frakS_4\cdot 2$. It leaves invariant 
the subset of $3$ points from $\Sing(\frakG)_1$ and $4$ points from $\Sing(\frakG)_2$. For example, we may assume 
that the plane is $\Pi = V(x_1+iy_1,x_2+iy_2,x_3-iy_3)$ and the points are 
$$[0,0,i,0,0,1],\  [0, -i,0,0,1,0],\ [-i, 0,0,1,0,0]$$
from the set $\Sing(\frakG)_1$ and 
$$[i, i, i, -1, -1, 1],\ [i, -i, i, -1, 1, 1], \ [i, i, -i, -1, -1, -1], \ 
[i, -i, -i, -1, 1, -1]$$
from the set $\Sing(\frakG)_2$. In the coordinates $y_1,y_2,y_3$, the seven points form the set of vertices of a complete quadrangle. 
Its group of symmetries is isomorphic to $\frakS_4$. 
On the other hand, a direct calculation shows that the pointwise stabilizer subgroup of
$\Pi$ is isomorphic to $\bbZ/2\bbZ$ generated by
$(x_1,x_2,x_3,y_1,y_2,y_3)\mapsto (-y_1,-y_2,y_3,x_1,x_2,-x_3)$ a composition of $g_0$ and
a sign change.  
This proves the assertion.
\end{proof}

Here we encounter the relationship with the $F_4$ root system again.  
As we explained in Section 3, the group of projective symmetries of the desmic pencil is isomorphic 
to a subgroup $W(F_4)/(w_0) \cong \half (\frakS_4\times \frakS_4)\rtimes 2$ of index 2 of $(\frakS_4\times \frakS_4)\rtimes 2$.

\vskip5pt
Consider the projection from $P=[0,0,0,0,0,1]$.  By putting $z_1 = p_{12}, z_2=p_{13},
z_3 = p_{14}, z_4 = p_{23}, z_5 = p_{24}, z_6 = p_{34}$ and substituting
$z_1z_6 = z_2z_5 - z_3z_4$ into (\ref{eqCubic3}), we obtain the quartic 3-fold
\beq
Y:z_1^2z_2z_4 -z_1^2z_3z_5 + z_2^2z_3z_5 - z_2z_3^2z_4 + z_3z_4^2z_5 - z_2z_4z_5^2 =0.
\eeq
We can rewrite the equation in the form
$$z_1^2(z_2z_4-z_3z_5)+(z_2z_3-z_4z_5)(z_2z_5-z_3z_4) = 0.$$
Let 
$$Q_1 = V(z_2z_4-z_3z_5), \quad Q_2 = V(z_2z_3-z_4z_5), \quad Q_3 = V(z_2z_5-z_3z_4).$$
The projection from the node $q = [1,0,0,0,0]$ of $Y$ defines a degree two map
$$\phi: \Bl_q(Y) \to V(z_1) \cong \bbP^3$$
branched along the union $Q_1\cup Q_2\cup Q_3$. 

The union of the quadrics $Q_2\cup Q_3$ is equal to the hyperplane section $V(z_1)\cap Y$.
They intersect along a quadrangle of lines  
$$V(z_1,z_2,z_4), \ V(z_1,z_3,z_5), \ V(z_1,z_2-z_4,z_3-z_5), \ V(z_1,z_2+z_4,z_3+z_5), $$
which are singular lines on $Y$. 
The intersections $Q_1\cap Q_2$ and $Q_1\cap Q_3$ are also equal to the union of four lines.

The quartic threefold $Y$ also has $17$ isolated ordinary nodes lying outside of the hyperplane 
$V(x_1)$.
\beq
\scalebox{0.7}{$\begin{split}
&[1,0,0,0,0], \ [1,0,0,1,1],\ [1,0,0,-1,1], \ [1,1,0,0,1], \ [1,-1,0,0,1],\\
&[1,0,0,1,-1], \ [1,0,0,-1,-1], \ [1,1,0,0,-1], \ [1,-1,0,0,-1], \ [1,0,1,1,0]\\
&[1,0,1,-1,0], \ [1,1,1,0,0], \ [1,-1,1,0,0], \ [1,0,-1,1,0], \ [1,0,-1,-1,0]\\
&[1,1,-1,0,0], \ [1,-1,-1,0,0].
\end{split}$}
\eeq
Note that sixteen of the nodes on $\frakG$ are projected to the points on four singular lines of $Y$. 
That explains the number 
 $17 = 33-16$. Under the projection $\phi$ from the first node $[1,0,0,0,0]$, the remaining sixteen nodes are 
 paired into eight pairs, and each pair is projected by one of the singular points  
 of $Q_1\cap Q_3$ or $Q_1\cap Q_3$.
 
\begin{theorem}\label{rationality}
The cubic complex $\frakG$ is a rational $3$-fold.
\end{theorem}
\begin{proof} 
It suffices to prove that the quartic 3-fold $Y$ is rational.
The proof is similar to that of the Burkhardt quartic three-fold \cite[\S 6]{Baker} (or see \cite[\S 5.2.7]{Hunt}).

Consider the following three planes contained in $Y$:
 \[
 \Pi_1 = V(z_2, z_3), \ \Pi_2=V(z_4,z_5), \ \Pi_3=V(z_1-z_2+z_4, z_1-z_3+z_5).
 \]
 Note that 
 \[
 \Pi_1\cap \Pi_2 = [1,0,0,0,0],\ \Pi_2\cap \Pi_3= [1,1,1,0,0],\ \Pi_3\cap \Pi_1=[1,0,0,-1,-1].
 \]
Let $P \in Y$ be a general point.  Then there exists a {\it unique}
 line $L$ meeting three planes $\Pi_1, \Pi_2,\Pi_3$ and passing through $P$, that is,
\[L = \la P, \Pi_1\ra \cap \la P, \Pi_2\ra \cap \la P,\Pi_3\ra,\]
where $\la P, \Pi_i\ra$ is the hyperplane generated by $P$ and $\Pi_i$.
Let $H$ be a hyperplane in $\bbP^4$.  Then we have a rational map
from $Y$ to $H$ by sending $P$ to $L\cap H$, which is obviously birational.
\end{proof}

 \begin{remark} A general hyperplane section of the Humbert desmic cubic line complex 
 is a K3 surface of degree $6$ in $\bbP^4$. It contains two sets of twelve disjoint 
 lines forming the configuration $(12_6,12_6)$. The geometry of these surfaces is 
 studied in great detail in  \cite{DDK}.
 \end{remark}

A complete intersection $V_{2,3}$ of a quadric $V(f_2)$ and a cubic $V(f_3)$ in $\bbP^n$ defines a cubic hypersurface 
$X_3$
in $\bbP^{n+1}$ given by the equation
\beq\label{cubic1}
x_0f_2(x_1,\ldots,x_{n+1})+f_3(x_1,\ldots,x_{n+1}) = 0.
\eeq
The point $[1,0,\ldots,0]$ is a double point $P$ of $X_3$ with the tangent cone $V(f_2)$. The closure of the union of lines 
through $P$ is a cone over $V(f_2,f_3)$. 
Conversely,  a cubic hypersurface $X_3$ with a double point $P$ can be given by equation \eqref{cubic1};
the variety $V(f_2,f_3)$ is called the \emph{associated variety} of 
$X_3$ and denoted by $\mathrm{AV}(X_3,P)$.

In general, a hypersurface singularity is an ordinary node ($A_1$-singularity) if and only if its corank is 0.  In our situation, the corank of any singular point on 
$\mathrm{AV}(X_3,P)$ coincides with
one of the corresponding singularities on $X_3$.  
In particular, if $X_3$ has only ordinary nodes, their number 
is equal to the number of ordinary nodes of $\mathrm{AV}(X_3,P)$ 
plus one (Kalker \cite[Cor. 1.15]{Kalker}).

Applying this to our situation, where $V_{2,3}$ is 
the Humbert cubic line complex $\frakG$ given by equation \eqref{humbertcomplex},
we obtain a cubic hypersurface in $\bbP^6$ given by the equation
\beq\label{ourcubic}
x_0(x_1^2+x_2^3+x_3^2+x_4^2+x_5^2+x_6^2)+x_1x_2x_3+\sqrt{-1}x_4x_5x_6 = 0
\eeq
which has 35 ordinary nodes.

 It is known that the maximum possible number of ordinary nodes on a nodal cubic hypersurface in $\bbP^n$ is equal to 
$\mu_n(3) = \binom{n+1}{[\frac{n}{2}]}$ \cite{Goryunov}, \cite{Kalker}. 
 The bound is achieved 
for the Segre cubic hypersurfaces
\beq\label{segrecubic}
t_0^3+\cdots+t_{2k+1}^3 = t_0+\cdots+t_{2k+1} = 0,
\eeq
if $n = 2k$, and the Clebsch-Segre cubic hypersurfaces
\beq\label{kalker}
t_0^3+\cdots+t_{2k+1}^3+2t_{2k+2}^3 = t_0+\cdots+t_{2k+1}+2t_{2k+2} = 0,
\eeq
if $n = 2k+1$ (it is projectively equivalent to the Goryunov cubic from \cite[Proposition 15]{Goryunov}).

In our case, where $n = 6$, 
we get $\mu_6(3) = 35$. 
This proves the following theorem.

\begin{theorem}\label{maximum}
We have the following assertions.
\begin{itemize}
\item[(1)] The number $34$ of ordinary nodes on the Humbert cubic line complex is 
the maximum possible number of ordinary nodes on a complete intersection $V_{2,3}$ 
in $\bbP^5$.  Moreover, the Humbert cubic line complex is the associated variety 
with respect to any node of
the Segre cubic hypersurface in $\bbP^6$.
\item[(2)] The Humbert cubic complex $\frakG$ of lines in $\bbP^3$ is projectively isomorphic to the associated variety 
of the Segre cubic hypersurface in $\bbP^6$.
\end{itemize}
\end{theorem}
\begin{proof}
We have proved the first assertion of (1) above, and the second one follows from 
the fact that the group $\mathfrak{S}_8$ of automorphisms of the Segre cubic 
hypersurface acts transitively on the set of its nodes.
To prove the second assertion (2), it is enough to show that the cubic hypersurface given by the equation 
\eqref{ourcubic} is projectively isomorphic to the Segre cubic hypersurface given by the equation \eqref{segrecubic}, 
where $k= 3$.
 
 The change of variables is given by 
\[
{\small \begin{split}
& t_0= 2x_0+x_1-x_2-x_3, \ t_1 = 2x_0-x_1+x_2-x_3,\ t_2 = 2x_0-x_1-x_2+x_3, \\
& 
t_3 = 2x_0+x_1+x_2+x_3,\ t_4 = -2x_0-i(x_4+x_5+x_6), \\
& t_5 = -2x_0+i(x_4+x_5-x_6),
\  t_6 = -2x_0+i(x_4-x_5+x_6),\\  
&t_7 = -2x_0+i(-x_4+x_5+x_6).
\end{split}}
\]
Thus, we have proved the assertion.
\end{proof}

\section{Cubic surfaces and desmic quartic surfaces}
Let $F$ be a smooth cubic surface and let $\ell_1+\ell_2+\ell_3$ be one of the triangles of lines on $F$ 
cut out by one of 
the $45$ tritangent planes. 

\begin{proposition} {\rm (L. Cremona)} Let $K_1,K_2,K_3$ be the residual conics for planes 
$\Pi_1,\Pi_2,\Pi_3$ from the
 pencils of planes through the lines $\ell_1,\ell_2,\ell_3$, respectively.
Then $K_1+K_2+K_3$ is cut out by a quadric $G(\Pi_1,\Pi_2,\Pi_3)$. The locus of singular points of 
quadrics $G(\Pi_1,\Pi_2,\Pi_3)$
is a quartic surface in $\bbP^3$.
\end{proposition}

\begin{proof}

Choose a geometric basis $(e_0,e_1,\ldots,e_6)$ in $\Pic(F)$ such that 
$[\ell_1] = e_0-e_1-e_2,\  [\ell_2]= e_0-e_3-e_4,\  [\ell_3] = e_0-e_5-e_6$. Then, 
the conics belong to $2e_0-(e_1+\cdots+e_6)+e_1+e_2$, etc. This implies 
that the sum is equal to $6e_0-2(e_1+\cdots+e_6) = -2K_F$. This proves the first assertion. 

Take a general line $\ell$ in $\bbP^3$. For each point $x\in \ell$, there exists a 
unique plane $\Pi_1(x),\Pi_2(x),\Pi_3(x)$ from each pencil of
 planes that 
intersects $\ell$ at $x$. This allows us to identify the three pencils of planes with the line $\ell$. Consider the map
$$\bbP^1\times \bbP^1\times \bbP^1\to |\calO_{\bbP^3}(2)|, \quad (\Pi_1,\Pi_2,\Pi_3)\mapsto G(\Pi_1,\Pi_2,\Pi_3).$$
Obviously, the map factors through $(\bbP^1\times \bbP^1\times \bbP^1)/\frakS_3\cong \bbP^3$ and defines a linear map
$$f:\bbP^3\to |\calO_{\bbP^3}(2)|.$$
Its image is a web of quadrics in $\bbP^3$. It is known that its Steinerian surface is a quartic surface 
\cite[Section 1.1.6]{CAG}. This proves the second assertion.
\end{proof}

We call the family of quadrics from the proposition the \emph{Cremona family of quadrics} associated with a tritangent plane.

Explicitly, fix a tritangent plane $w=0$ intersecting $F$ along the lines $\ell_1 = V(x,w), \ell_2 = V(y,w), \ell_3 = V(z,w)$. 
Here we assume that the tritangent plane is non-degenerate, that is, three lines do not meet
in a point (however, we have the same equation \eqref{eqnquadric} for the degenerate case). 
The equation of $F$ can be written in the form
\beq\label{eqncubic}
q(x,y,z,w)w+xyz= 0,
\eeq
 where $q$ is a quadratic form. 
 The pencil of planes through the line $\ell_1$, $\ell_2$ or $\ell_3$ 
is given by $V(x-\alpha w)$, $V(y-\beta w)$ or $V(z-\gamma w)$, respectively. 
Plugging in the equation of $F$, we find that the conics are 
$V(q+\alpha yz,x-\alpha w), V(q+\beta xz,y-\beta w)$, and $V(q+\gamma xy,z-\gamma w)$. We check 
that the conics lie in the quadric
\beq\label{eqnquadric}
q+\alpha yz+\beta xz+\gamma xy-(\alpha\beta z+\alpha\gamma y+\beta\gamma x)w+\alpha\beta\gamma w^2 = 0.
\eeq
By homogenizing, we see that its coefficients are tri-linear functions of 
$(\alpha_0,\alpha_1)$, $(\beta_0,\beta_1)$, $(\gamma_0,\gamma_1)$. 
The locus of singular points of these quadrics is given by the equations
\begin{equation}\label{eqnjacobian}
\begin{split}
&q_x+\beta z+\gamma y-\beta \gamma w = 0,\\
&q_y+\alpha z+\gamma x-\alpha \gamma w = 0,\\
&q_z+\beta x+\alpha y-\alpha \beta w = 0,\\
&q_w-\alpha\beta z-\alpha\gamma y-\beta\gamma x+2\alpha\beta\gamma w =0.
\end{split}
\end{equation}
 The following equation of the Steinerian surface was found by Humbert:
\beq\label{humbert}
G = q^2 -q_yq_zyz-q_xq_zxz-q_xq_yxy-q_xq_yq_zw+xyzq_w = 0.
\eeq
It can be checked by using the equations \eqref{eqnquadric} and \eqref{eqnjacobian}.
 
We can find a normal form for a cubic surface \eqref{eqncubic} with a fixed tritangent plane. 
The only admissible transformations 
are $(x,y,z,w)\mapsto (\lambda x+c_1w,\mu y+c_2w,\gamma z+c_3w,\delta w)$, where $\lambda\mu\gamma = 1$. 
Using these transformations, we may assume that 
\beq\label{eqncubic2}
((aw + bx + cy + dz)w + x^2 + y^2 + z^2)w+xyz= 0.
\eeq 
This agrees with the fact that cubic surfaces depend on four parameters. This gives the equation of Cremona's 
quartic surface
\beq\label{eqncremona}
\begin{split}
&-(bw + 2x)(cw + 2y)(dw + 2z)w - (bw + 2x)(cw + 2y)xy\\
 &- (bw + 2x)(dw + 2z)xz
  - (cw + 2y)(dw + 2z)yz  
+(2aw + bx + cy + dz)xyz \\
&+ (aw^2 + bwx + cwy + dwz + x^2 + y^2 + z^2)^2 =0.
\end{split}
\eeq
Substituting the expression for $q_x,q_y,q_z,q_w$ in equations \eqref{eqnjacobian}, we find the equation 
of the discriminant surface of the family of quadrics \eqref{eqnquadric}:
\beq
\det\left(\begin{smallmatrix}2\alpha_0\beta_0\gamma_0&\alpha_0\beta_0\gamma_1&\alpha_0\beta_1\gamma_0
&-\alpha_0\beta_1\gamma_1+b\alpha_0\beta_0\gamma_0\\
\alpha_0\beta_0\gamma_1&2\alpha_0\beta_0\gamma_0&\alpha_1\beta_0\gamma_0&-\alpha_1\beta_0\gamma_1+c\alpha_0\beta_0\gamma_0\\
\alpha_0\beta_1\gamma_0&\alpha_1\beta_0\gamma_0&2\alpha_0\beta_0\gamma_0&-\alpha_1\beta_1\gamma_0+d\alpha_0\beta_0\gamma_0\\
-\alpha_0\beta_1\gamma_1+b\alpha_0\beta_0\gamma_0&-\alpha_1\beta_0\gamma_1+c\alpha_0\beta_0\gamma_0
&-\alpha_1\beta_1\gamma_0+d\alpha_0\beta_0\gamma_0&2(a\alpha_0\beta_0\gamma_0+\alpha_1\beta_1\gamma_1)\end{smallmatrix}\right) = 0.
\eeq
Here, we homogenized the affine coordinates $\alpha,\beta,\gamma$ in $\bbP^1\times \bbP^1\times \bbP^1$.
The determinantal equation is given by a multi-homogeneous form of multi-degree $(4,4,4)$, which is symmetric with respect 
to permuting the factors of $(\bbP^1)^3$. After we identify $(\bbP^1)^3/\frakS_3$ with $\bbP^3$, we obtain a symmetric 
determinant whose entries are linear functions in homogeneous coordinates. This is the discriminant of Cremona's 
web of quadrics.  Its Steinerian surface is given by \eqref{eqncremona}.

The following theorem is due to Cremona \cite[p. 91]{Cremona} and Steiner \cite[p. 80]{Steiner}. 
We give a proof following Humbert \cite[Part IV]{Humbert}.

\begin{theorem}\label{12-4,16-3} The Steinerian quartic surface $S$ has $12$ nodes 
and $16$ lines forming a Reye configuration $(12_4,16_3)$.
\end{theorem}

\begin{proof} Fix a tritangent plane $\Pi_0 = V(w)$ as above. 
Through each line $\ell_i$ passes $5$ tritangent planes including $\Pi_0$.
This gives $24 = 12\times 2$ lines in these planes not counting $\ell_1,\ell_2,\ell_3$. There are $45-12-1 = 32$ other tritangent planes.
They are distributed in 16 pairs. Each pair contains six of the lines from the set of 24 planes and defines a reducible quadric from 
the Cremona family of quadrics \cite[p. 91]{Cremona}.   

For example, in the Cremona's notation 
\[
a_i = e_i,\quad b_j = 2e_0-e_1 -\cdots - e_6 + e_j, \quad c_{ij}=e_0 -e_i-e_j
\] 
of 27 lines on a cubic surface (see \cite[Section 9.1.1]{CAG}), take the tritangent plane $T$ with lines 
$\{c_{12},c_{34},c_{56}\}$. Consider three tritangent planes $T_1,T_2,T_3$ with lines 
$\{c_{12},c_{35},c_{46}\}$, $\{c_{16},c_{25},c_{34}\}$, $\{c_{14},c_{23},c_{56}\}$, respectively. 
Each contains one line from $T$. The pair of tritangent planes $T_4,T_5$ with lines 
$\{c_{13},c_{26},c_{45}\}$, $\{c_{15},c_{24},c_{36}\}$ intersect each of the planes $T_1,T_2,T_3$ along a line. The union of these tritangent planes is a reducible 
quadric that cuts out in each plane $T_1, T_2, T_3$ a reducible conic. Note that the tritangent planes 
$T, T_4, T_5$ form one of the conjugate triads of tritangent planes (see \cite[Section 9.1.1]{CAG})
$$\begin{matrix}c_{12}&c_{34}&c_{56}\\
c_{36}&c_{15}&c_{24}\\
c_{45}&c_{26}&c_{13}.\end{matrix}$$
The tritangent planes $T_1, T_2, T_3$ are formed by the remaining $15-9 = 6$ lines $c_{ij}$. 
Similarly, we find other conjugate triads with 
the top row $c_{12},c_{34},c_{56}$, for example
$$\begin{matrix} c_{12}&c_{34}&c_{56}\\
a_1&b_3&c_{13}\\
b_2&a_4&c_{24}.\end{matrix}
$$
The tritangent planes defined by the last two rows define a reducible quadric that cuts out line-pairs in the 
tritangent planes with lines $\{c_{12},a_3,b_1\}$, $\{c_{34},a_3,b_4\}$, $\{c_{56},c_{14},c_{23}\}$.
Proceeding in this way, we find 16 conjugate triads of tritangent planes with the top row 
$\{c_{12},c_{34},c_{56}\}$ that give us 16 reducible quadrics from the Cremona family. The singular line of each reducible quadric 
is a line on the Steinerian surface.

To find 12 double points, we use \eqref{humbert}. Let $f:=qw+xyz = 0$ be the equation of $F$. We find 
$$f^2-f_xf_yf_z = Gw^2.$$
The tangent plane to $F$ at a point $[x_0,y_0,z_0,w_0]$ is given by equation 
\[
f_x(x_0,\ldots,w_0)x+\cdots+f_w(x_0,\ldots,w_0)w= 0.
\]
This implies that the intersection points $V(f,f_x,f_y)$ are the points where the lines $V(t_0z+t_1w)$ are tangent 
to $F$.  Each of the sixteen lines from above intersects the cubic surface at three points. 
The points are the intersection points 
with the three lines lying in each irreducible component of the plane-pair intersecting along this line. 
This irreducible component 
is a tritangent plane containing one of the lines $\ell_1,\ell_2,\ell_3$, for example, the line $V(z,w)$. 
It follows that these points belong to $V(f,f_x,f_y)$, and hence, they are double points of the Steinerian quartic $Q$.
The fact that the configuration of points and lines is isomorphic to the Reye configuration follows from the next theorem.
\end{proof}

The following theorem is due to Humbert. We give it another proof.

\begin{theorem} The quartic surface $Q$ associated with a tritangent plane of a smooth cubic surface $F$ is 
isomorphic to a desmic quartic surface.
\end{theorem}

\begin{proof} We found the same configuration $(12_4,16_3)$ of nodes and lines on the Steinerian quartic. 
Now, we follow the proof 
of Theorem \ref{mainthm1} to show that it is isomorphic to the Kummer surface.
\end{proof}

Here are the immediate questions:
\begin{itemize}
\item A cubic surface admits 45 tritangent planes. Are the corresponding desmic quartic surfaces isomorphic?
\item The construction assigns to a cubic surface and a tritangent plane $(F,T)$ on it the isomorphism class of an elliptic 
curve. Is there a more direct way to see this association?
\item What are the fibers of the map from the moduli space of cubic surfaces with a marked 
tritangent plane to the moduli space of elliptic curves? It is known that the former moduli space is isomorphic to the quotient of 
the blow-up of a maximal torus of a simple Lie group of type $F_4$ at the identity element by the Weyl group $W(F_4)$ (see
\cite{Looijenga}). Is there any relation to the $W(F_4)$ action on the desmic pencil and the cubic line complex?
\item What happens with the Cremona's quartic surface if the tritangent plane contains an Eckardt point?

\end{itemize}

\section{Cremona's quartic surfaces in characteristic two}
Starting from this section, we assume that the characteristic of $\Bbbk$ is equal to $2$.
Let us rewrite the equation of Cremona's quartic surface from \eqref{eqncremona} in characteristic 2. The equation 
of the cubic surface \eqref{eqncubic} still holds. The discriminant surface 
of the web of quadrics is given now by the pfaffian of the alternating matrix:
\beq
\textrm{Pf}\left(\begin{smallmatrix}0&\alpha_0\beta_0\gamma_1&\alpha_0\beta_1\gamma_0
&\alpha_0\beta_1\gamma_1+b\alpha_0\beta_0\gamma_0\\
\alpha_0\beta_0\gamma_1&&\alpha_1\beta_0\gamma_0&\alpha_1\beta_0\gamma_1+c\alpha_0\beta_0\gamma_0\\
\alpha_0\beta_1\gamma_0&\alpha_1\beta_0\gamma_0&0&\alpha_1\beta_1\gamma_0+d\alpha_0\beta_0\gamma_0\\
\alpha_0\beta_1\gamma_1+b\alpha_0\beta_0\gamma_0&\alpha_1\beta_0\gamma_1+c\alpha_0\beta_0\gamma_0
&\alpha_1\beta_1\gamma_0+d\alpha_0\beta_0\gamma_0&0\end{smallmatrix}\right) = 0.
\eeq
It is a multi-linear form of multi-degree $(2,2,2)$. So, it becomes a quadric in $(\bbP^1)^3/\frakS_3 = \bbP^3$.
The Humbert equation of Cremona's quartic still holds. It expresses the condition that 
the partial derivatives of a quadric from the web vanish, and also the quadric vanishes at the same point.
From (\ref{eqncremona}), we obtain the equation
\beq\label{eqncremona2}
\begin{split}
&G= bcdw^4 + bcw^2xy
 + bdw^2xz
  + cdw^2yz  
+(bx + cy + dz)xyz \\
&+ (aw^2 + bwx + cwy + dwz + x^2 + y^2 + z^2)^2 =0.
\end{split}
\eeq
The proof that the surface $V(G)$ contains 16 lines does not depend on the characteristic. Let us find the singularities.
We have
\beq
G_x' = (cy+dz)(yz+bw^2),\ G_y' = (bx+dz)(xz+cw^2), \ G_z'= (bx+cy)(xy+dw^2),
\eeq
and $G_w'=0$.
One checks that, for a general parameter, no singular points lie in the plane $w= 0$. So, we may 
assume that $w = 1$ and find the following singular points. Four points
$$[dz_1^2,b^2,bz_1^2,bz_1], \quad b^2dz_1^3 + b^2z_1^4 + d^2z_1^4 + ab^2z_1^2 + b^3cz_1 + b^4 = 0.$$ 
Four points 
$$[c^2,dz_2,cz_2^2,cz_2], \quad c^2dz_2^3 + c^2z_2^4 + d^2z_2^4 + ac^2z_2^2 + bc^3z_2 + c^4 = 0.$$ 
Four points
$$[c,b,z_3^2,z_3], \quad dz_3^3 + z_3^4 + az_3^2 + bcz_3 + b^2 + c^2 = 0.$$
There is also an additional singular  point in characteristic 2:
$$P_0 = [cdz_0,bdz_0,bcz_0,bc], \quad z_0^4 = \frac{(b^5c^5d + a^2b^4c^4)}{(b^3c^3d^3 + 
b^4c^4 + b^4d^4 + c^4d^4)}.$$
The line connecting a point from the first group of four singular points with 
a point from the second group also contains one of the singular points from the 
third group, and belongs to the surface.
This defines a Reye configuration $(12_4,16_3)$ formed by the singular points and sixteen lines.
So, we may keep the name desmic surface for Cremona's quartic $V(G)$, although no desmic pencil exists in characteristic two.
 Note that a minimal nonsingular model of the surface contains 16 disjoint 
 $(-2)$-curves, and hence, it is a supersingular K3 surface
 \cite{SB}. 
 
 \vskip5pt
Let us determine the type of the singular point $P_0$. For special parameters $(a,b,c,d) = (0,0,1,1)$ for 
which the cubic surface is smooth, we checked directly that $P_0$ is a rational double point of type $A_3$.

Let us show that this is always true if we assume that $P_0$ does not lie on any of the sixteen lines on the surface.
We will show later that a K3 surface containing a $(12_4,16_3)$-configuration is a supersingular K3 surface $X$ with Artin invariant $\sigma\leq 2$ (see Theorem \ref{Artin2}).
In this case, the lattice $L$ generated by twenty-eight $(-2)$-curves forming $(12_4,16_3)$-configuration can be primitively embedded in the Picard lattice $S_\sigma=\Pic(X)$.
Recall that $L=U\oplus D_8\oplus D_9\cong U\oplus D_5\oplus D_{12}$ (Proposition \ref{PicardLattice}).
In case $\sigma=1$, $S_\sigma = U\oplus E_8\oplus D_{12}$ is an over lattice of the
orthogonal direct sum of $L=U\oplus D_5 \oplus D_{12}$ and $A_3$ ($E_8$ is an overlattice of $D_5\oplus A_3$).  In case $\sigma=2$, $S_\sigma = U\oplus D_8\oplus D_{12}$ is an over lattice of the orthogonal direct sum of $L=U\oplus D_5 \oplus D_{12}$ and $A_3$ ($D_8$ is an overlattice of $D_5\oplus A_3$ which is the intermediate one: 
$D_5\oplus A_3\subset D_8\subset E_8$).  The linear system (\ref{pola}) gives a
quartic model which contracts twelve disjoint $(-2)$-curves and 
three $(-2)$-curves forming a root lattice of type $A_3$.  
Thus, the additional singularity is of type $A_3$.

\section{Kummer surfaces $\Kum(E\times E)$ in characteristic two}

In this section, we discuss what happens for Kummer surfaces in characteristic two, 
and in the next section, we will locate the Reye configuration $(12_4,16_3)$ 
among $42$ smooth rational curves on the 
supersingular K3 surface with Artin invariant $\sigma = 1$ studied in \cite{DK}.

Let $E$ be an ordinary elliptic curve in characteristic two and $\iota:E\to E$ 
the inversion map.
Denote by $\{0, a\}$ the 2-torsion points of $E$.  
Let
\[
E_1= E \times \{0\},\ E_2 =E \times \{a\},\ F_1 = \{0\}\times E,\ F_2 = \{a\} \times E,
\]
and let $\Delta_1$ be the diagonal and $\Delta_2$ the translation of $\Delta_1$ by a 2-torsion $(0,a)$.  The quotient surface $\Kum(E\times E)$ 
of $E\times E$ by the inversion involution $\iota\times \iota$ has 
four singular points of type $D_4$ \cite{Shioda}.  
Its minimal resolution $\widetilde{\Kum}(E\times E)$ has the following
22 $(-2)$-curves as in Figure \ref{ExEdualgraph}:

\begin{figure}[h]
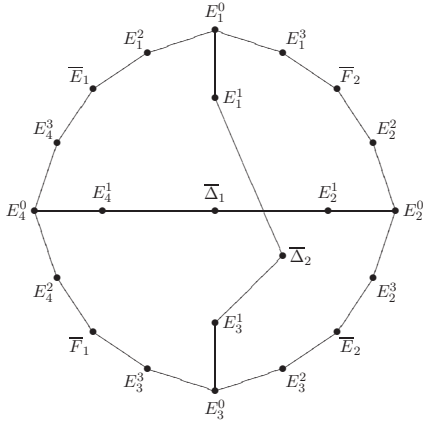

\begin{center}
\scalebox{0.6}{
\xy 
@={(40,0),(0,-40),(-40,0),(0,40),(27,27),(27,-27),(-27,27),(-27,-27),(15,35),(15,-35),(-15,35),
(-15,-35),(35,15),(35,-15),(-35,15),(-35,-15),(0,25),(0,-25),(-25,0),(25,0),(15,-10),(0,0)}@@{*{\bullet}};
(40,0)*{};(35,15)*{}**\dir{-};(35,15)*{};(27,27)*{}**\dir{-};(27,27)*{};(15,35)*{}**\dir{-};
(15,35)*{},(0,40)*{}**\dir{-};
(0,40)*{};(-15,35)*{}**\dir{-};(-15,35)*{};(-27,27)*{}**\dir{-};
(-27,27)*{};(-35,15)*{}**\dir{-};(-35,15)*{};(-40,0)*{}**\dir{-};
(-40,0)*{};(-35,-15)*{}**\dir{-};(-35,-15)*{};(-27,-27)*{}**\dir{-};
(-27,-27)*{};(-15,-35)*{}**\dir{-};
(-15,-35)*{},(0,-40)*{}**\dir{-};
(0,-40)*{};(15,-35)*{}**\dir{-};(15,-35)*{};(27,-27)*{}**\dir{-};
(27,-27)*{};(35,-15)*{}**\dir{-};(35,-15)*{};(40,0)*{}**\dir{-};
(40,0)*{};(25,0)*{}**\dir{-};(-40,0)*{};(-25,0)*{}**\dir{-};
(0,40)*{};(0,25)*{}**\dir{-};(0,-40)*{};(0,-25)*{}**\dir{-};
(25,0)*{};(0,0)*{}**\dir{-};
(-25,0)*{};(0,0)*{}**\dir{-};
(0,25)*{};(15,-10)*{}**\dir{-};
(0,-25)*{};(15,-10)*{}**\dir{-};
(44,0)*{E_2^0};(-44,0)*{E_4^0};(0,-44)*{E_3^0};(0,44)*{E_1^0};
(30,30)*{\overline{F}_2};(30,-30)*{\overline{E}_2};(-30,30)*{\overline{E}_1};(-30,-30)*{\overline{F}_1};
(18,38)*{E_1^3};(18,-38)*{E_3^2};(-18,38)*{E_1^2};(-18,-38)*{E_3^3};
(38,18)*{E_2^2};(38,-18)*{E_2^3};(-38,18)*{E_4^3};(-38,-18)*{E_4^2};
(4,25)*{E_1^1};(4,-25)*{E_3^1};(25,4)*{E_2^1};(-25,4)*{E_4^1};
(0,4)*{\overline{\Delta}_1};(19,-10)*{\overline{\Delta}_2};
\endxy
}
\end{center}
\caption{The dual graph: ordinary case}
\label{ExEdualgraph}
\end{figure}

Here, $\overline{E}_i, \overline{F}_i, \overline{\Delta}_i$ are the images of 
$E_i, F_i, \Delta_i$ and
$\{E_i^j\}_{j=0,1,2,3}$ are the exceptional curves over a singular point of type $D_4$ 
$(i=1,2,3,4)$.

Recall that a desmic surface in characteristic $p\ne 2$ is obtained from 
$\widetilde{\Kum}(E\times E)$ by contracting twelve $(-2)$-curves (Remark \ref{from Kummer}).
Let us consider its analog in characteristic two.
Let $p_1, p_2$ be the first and second projections and let
$p_3: E\times E \to E$, $(x,y)\to (x+y, x+y)$.
Each $p_i$ induces an elliptic fibration $\pi_i: \widetilde{\Kum}(E\times E)\to \bbP^1$ which
has two singular fibers of type $\tilde{D}_8$ (type ${\rm I}^*_4$ in Kodaira's notation). 
We denote by $\calF_i$ a general fiber of $\pi_i$.

Now, consider the divisor
\[
H = \calF_1+\calF_2+\calF_3 - \sum_{i=1}^4 \sum_{j=1}^3 E_i^j - 2\sum_{i=1}^4 E_i^0.
\]
It is easy to see that 
$H^2=4, H\cdot E_i^0 =1\ (i=1,2,3,4)$ and
other curves in Figure \ref{ExEdualgraph} are orthogonal to $H$.  
Thus, the linear system $|H|$ gives
a birational morphism from $\widetilde{\Kum}(E\times E)$ to a quartic surface $S$ with six rational double
points of type $A_3$ (= the images of 18 $(-2)$-curves in Figure \ref{ExEdualgraph} except $E_i^0$). 
Let $\overline{E}_i^0$ be the images of $E_i^0$, all of which are lines contained in a plane, 
and their six intersection points are rational double points of $S$ of type $A_3$.

Conversely, let 
\[
F_{\alpha} = xyzw + \alpha (x+y+z+w)^4,\ \alpha\ne 0 \in k.
\]
Then, $Q_{\alpha}=V(F_{\alpha})$ contains four lines, irreducible components of $V(xyzw)$ 
in the plane $V(x+y+z+w)$ and has six singular points at the six intersection points of four lines:
\[
[1,1,0,0],\ [1,0,1,0],\ [1,0,0,1],\ [0,1,1,0],\ [0,1,0,1],\ [0,0,1,1].
\]
Since each singular point is formally analytically given by $t^4=uv$, it is a rational double point of type $A_3$.
Thus, the four lines and six singular points form a $(6_2, 4_3)$-configuration.
Contrary to the case of characteristic not equal to 2, only two (degenerate) tetrahedra
$V(xyzw), V((x+y+z+w)^4)$ appear in the pencil $\{Q_{\alpha}\}_{\alpha \in\bbP^1}$.

Let $X_{\alpha}$ be the minimal resolution of singularities of $Q_{\alpha}$.
Then $X_{\alpha}$ contains twenty-eight $(-2)$-curves forming the dual graph 
in Figure \ref{ExEdualgraph}.
It follows from \cite[Theorem 3.4]{SB} that $X_{\alpha}$ is 
birationally isomorphic to a Kummer surface, that is,
there exists an abelian surface $A$ with $X_{\alpha}= \widetilde{\Kum}(A)$.  
The existence of three elliptic fibrations with two singular fibers of type $\tilde{D}_8$
as $\pi_1, \pi_2, \pi_3$ implies that $A\cong E\times E$ for an ordinary elliptic curve $E$.  
We now conclude the following theorem.

\begin{theorem}\label{KumChar2}
The quartic surface $Q_{\alpha}=V(F_{\alpha})$ is  birationally isomorphic to the Kummer surface $\Kum(E\times E)$
for an ordinary elliptic curve $E$.
\end{theorem}

\begin{remark}
A pencil of planes containing a line on $Q_\alpha$ induces an elliptic fibration on
$X_{\alpha}$ with singular fibers of type $\tilde{E}_6$ and of type $\tilde{A}_{11}$.
\end{remark}

\begin{remark}
The automorphism group of $\widetilde{\Kum}(E\times E)$ is calculated in \cite{KoMu}. It coincides with the group of birational automorphisms 
of a desmic quartic surface.
\end{remark}

\begin{figure}[htbp]
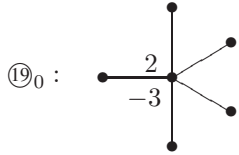

\begin{center}
\resizebox{!}{2.0cm}{
\xy
(80,15)*{};(90,15)*{}**\dir{-};
(90,15)*{};(90,25)*{}**\dir{-};
(90,15)*{};(90,5)*{}**\dir{-};
(90,15)*{};(98.5,20)*{}**\dir{-};
(90,15)*{};(98.5,10)*{}**\dir{-};
@={(80,15),(90,5),(90,15),(98.5,10),
(98.5,20),(90,25)
}@@{*{\bullet}};
(86,12)*{-3};
(87,17)*{2};
(70,15)*{\mbox{\MARU{19}}_{0}:};
\endxy 
}
\end{center}
\caption{The dual graph of a minimal resolution of singularity}
\label{19_0}
\end{figure}

In the case where $E$ is supersingular, the equation of the quotient $E\times E/\la \iota\times \iota\ra$ is given by
\[
z^2 +x^2y^2z + x^4y + xy^4 =0,
\]
and it has an elliptic singularity 
of type $\mbox{\MARU{19}}_0$ in the notation from \cite{Wagreich}
(see \cite{Katsura}).  The dual graph of its minimal resolution is given in Figure \ref{19_0}.
The central component has multiplicity 2 and self-intersection number $-3$.  
Other components are $(-2)$-curves.

We do not know how to find a family of quartic surfaces that contains the Kummer 
surfaces of the self-product of the supersingular elliptic curve.

\section{Supersingular K3 surfaces in characteristic two}

In this section, we first show that the supersingular K3 surface $X$ with Artin invariant
$1$ in characteristic two has twenty-eight $(-2)$-curves forming $(12_4,16_3)$-configuration explicitly.  Next, 
we prove that a supersingular K3 surface with
Artin invariant $\sigma$ in characteristic two has such $(-2)$-curves if and only if $\sigma \leq 2$.

First, let us recall that $X$ is obtained as the minimal resolution of an inseparable double covering $\bar{X}:
w^2 = xyz(x^3+y^3+z^3)$ of $\bbP^2$ \cite{DK}.  
The covering surface $\bar{X}$ has 21 nodes over 21 points from $\bbP^2(\bbF_4)\subset \bbP^2(\Bbbk)$.
Thus $X$ contains two families of $21$ disjoint $(-2)$-curves.
One of them consists of the exceptional curves and the other consists of proper transforms of 21 lines on $\bbP^2(\bbF_4)$. 
Each member of one family meets exactly five members in another family, that is,
they form a $(21_5)$-configuration.

Next, we recall Sylvester's duads and synthemes (see \cite[Section 9.4.3]{CAG}). 
First of all, we fix a set of six letters $\{1, 2, \cdots, 6\}$.
We denote by $ij$ the pair of $i$ and $j$ 
($1 \leq i \not= j \leq 6$) which is classically called {\it Sylvester's duad}.
Six letters $1,2,3,4,5,6$ can be arranged in three pairs of duads,
for example, $12.34.56$, called {\it Sylvester's syntheme}.  
(It is understood
that $12.34.56$ is the same as $12.56.34$ or $34.12.56$).
Duads and synthemes form a symmetric $(15_3)$-configuration.  
It is possible
to choose a set of five synthemes that together contain all the
fifteen duads.  Such a family is called a {\it total}.  The number of 
possible totals is six.  And every two totals have one, and only one
syntheme in common between them.  The following table gives 
the six totals in its row, and also in its columns (see \cite[p. 465]{CAG}) :
\bigskip
\bigskip

\halign{\hfil\tt#\hfil&&\quad#\hfil\cr
&       & $T_1$ & $T_2$ & $T_3$ & $T_4$ & $T_5$ & $T_6$ \cr
& $T_1$ &  & 14.25.36 & 16.24.35 & 13.26.45 & 12.34.56& 15.23.46 \cr
\noalign{\smallskip}
& $T_2$ & 14.25.36& & 15.26.34& 12.35.46& 16.23.45& 13.24.56\cr
\noalign{\smallskip}
& $T_3$ & 16.24.35& 15.26.34&  & 14.23.56& 13.25.46& 12.36.45\cr
\noalign{\smallskip}
& $T_4$ & 13.26.45& 12.35.46& 14.23.56& & 15.24.36& 16.25.34\cr
\noalign{\smallskip}
& $T_5$ & 12.34.56& 16.23.45& 13.25.46& 15.24.36& & 14.26.35\cr
\noalign{\smallskip}
& $T_6$ & 15.23.46& 13.24.56& 12.36.45& 16.25.34& 14.26.35& \cr}

\medskip

 We can label the sets of
21 points,  21 lines on $\bbP^2(\bbF_4)$, six totals $T_1,\cdots, T_6$,
fifteen duads, fifteen syntemes as follows (see \cite[Chap. 23, \S 4]{CS}).  We take six points on $\bbP^2(\bbF_4)$ in 
general position,
i.e., no three of them that are collinear, and identify them with $1,\cdots, 6$.  Then, we identify
the line passing through $i, j$ with the duad $ij$.  The remaining three points on the line $ij$, except
$i, j$ can be identified with three synthemes containing the duad $ij$.  Finally, we identify 
the remaining six lines with the six totals.
Thus, we can label 42 $(-2)$-curves on $X$, too.  
We remark that the incidence relation between lines and points corresponds to that between points, totals, duads and synthemes, that is, a point in $\bbP^2(\bbF_4)$ is contained in a line if and only if 
the corresponding point from $\{1,\cdots, 6\}$ (resp. or a syntheme) is contained in the corresponding
duad (resp. or a total).

Now, we show the existence of the same type of three genus one fibrations $f_1, f_2, f_3$ as
in Remark \ref{from Kummer}.  In this case, they are quasi-elliptic.
First, take the duad $12$.  The four $(-2)$-curves $12.35.46$, $12.34.56$, $12.36.45$, $2$ together with
$12$ form a singular fiber of type $\tilde{D}_4$ in which $12$ is the central irreducible component.
This cycle defines a genus one fibration $f_1$ with five singular fibers of type $\tilde{D}_4$.
This is quasi-elliptic by \cite[\S 1, Proposition]{RS}.
Similarly, we define $f_2$ (resp. $f_3$) by $26$ (resp. $23$).
In columns of Tables \ref{T1}, \ref{T2}, \ref{T3}, we give five singular fibers of $f_1, f_2, f_3$,
respectively.

\begin{table}[ht]
 \centering
  \begin{tabular}{|c|c|c|c||c|}
   \hline  
 $12$ & $13$ & $14$& $15$ & $16$ \\
 \hline \hline
  $12.35.46$ & $13.25.46$ & $4$ & $15.23.46$ & $6$\\ 
 $12.34.56$  & $13.24.56$  & $14.23.56$ & $5$ & $16.23.45$\\ 
 $12.36.45$ & $3$ & $14.25.36$ & $15.24.36$ & $16.24.35$\\ 
$2$ & $13.26.45$ &  $14.26.35$  & $15.26.34$ & $16.25.34$ \\
\hline
\end{tabular}
\caption{Singular fibers of $f_1$}
\label{T1}
\end{table}

\begin{table}[ht]
 \centering
  \begin{tabular}{|c|c|c|c||c|}
   \hline  
$26$ & $36$ & $46$& $56$ & $16$ \\
 \hline \hline
 $2$ & $12.36.45$ & $12.35.46$ & $12.34.56$ & $1$\\ 
$13.26.45$  & $3$  & $13.25.46$ & $13.24.56$ & $16.23.45$\\ 
 $14.26.35$ & $14.25.36$ & $4$ & $14.23.56$ & $16.24.35$\\ 
 $15.26.34$ & $15.24.36$ & $15.23.46$  & $5$ & $16.25.34$ \\
\hline
\end{tabular}
\caption{Singular fibers of $f_2$}
\label{T2}
\end{table}

\begin{table}[ht]
 \centering
  \begin{tabular}{|c|c|c|c||c|}
   \hline  
$23$ & $45$ & $T_2$& $T_5$ & $16$ \\
 \hline \hline
 $2$ & $4$ & $12.35.46$ & $12.34.56$ & $1$\\ 
$3$  & $5$  & $13.24.56$ & $13.25.46$ & $6$\\ 
 $14.23.56$ & $12.36.45$ & $14.25.36$ & $14.26.35$ & $16.24.35$\\ 
 $15.23.46$ & $13.26.45$ & $15.26.34$  & $15.24.36$ & $16.25.34$ \\
\hline
\end{tabular}
\caption{Singular fibers of $f_3$}
\label{T3}
\end{table}

Note that sixteen simple components of the first four singular fibers are common in $f_1$, $f_2$, and $f_3$.
We consider these sixteen disjoint $(-2)$-curves and 
twelve disjoint $(-2)$-curves $12$, $13$, $14$, $15$, $26$, $36$, $46$, 
$56$, $23$, $45$, $T_2$, $T_5$
(central components of the first four singular fibers of $f_1$, $f_2$, $f_3$).  
The incidence relation between these twenty-eight curves is the same as in Remark \ref{from Kummer}.
The linear system $|H|$ given in (\ref{pola}) defines a birational morphism from $X$ to a quartic surface
such that the twelve $(-2)$-curves are contracted to nodes, and the images of the sixteen curves are lines.
The $(-2)$-curves not occurring here are the following $14 (=42-28)$ curves:
\[
1,\ 6,\ 16,\ 24,\ 25,\ 34,\ 35, \ 16.23.45,\ 16.24.35,\ 16.25.34,\ T_1,\ T_3,\ T_4,\ T_6.
\]
By calculating their intersection numbers with $H$, we can see that the images of 
$1$, $6$, $16.23.45$ are conics, those of $24$, $25$, $34$, $35$, $T_1$, $T_3$, $T_4$, $T_6$ are cubics 
and $16$, $16.24.35$, $16.25.34$ are contracted to a rational double point of type $A_3$.
Thus, we have the following theorem.

\begin{theorem}\label{Artin1}
Let $X$ be the supersingular $K3$ surface with Artin invariant $\sigma =1$ 
in characteristic two.  Then, $X$  contains twenty-eight $(-2)$-curves forming 
$(12_4, 16_3)$-configuration.  The surface $X$ admits a quartic birational model, 
which has sixteen lines, twelve nodes, and a rational double point of type $A_3$.  
\end{theorem}

Next, we discuss when a supersingular K3 surface with Artin invariant $\sigma$ in characteristic two contains such
twenty-eight $(-2)$-curves.  
To do this, we use the quadratic lattice theory due to Nikulin \cite{Nikulin2}.  

Let $L$ be a lattice, i.e., $L$ is a free $\bbZ$-module of finite rank together with
a non-degenerate $\bbZ$-valued symmetric bilinear form $\la , \ra : L\times L\to \bbZ$. 
A lattice $L$ is called {\it even} if $\la x, x\ra$ is even for all $x\in L$.  
Since $L$ is non-degenerate, the dual 
$L^* = \Hom(L,\bbZ)$ of $L$ contains $L$.   Let $A_L=L^*/L$, 
which is a finite abelian group.  We denote by $l(L)$ the number of minimal generators of $A_L$. 
For an even lattice $L$, define 
$q_L: L^*/L \to \bbQ/2\bbZ$, by
$q_L(x+L) = \la x, x\ra +2\bbZ$, which is called the {\it discriminant quadratic form}.
Let $S$ be an even lattice and let $L\subset S$ be a {\it primitive} sublattice, 
i.e. $S/L$ is torsion free.  Let $M = L^\perp$ be the orthogonal complement of $L$, 
which is a primitive sublattice of $S$.  
Then $H_S = S/L\oplus M \subset A_L\oplus A_M$ is an isotropic
subgroup with respect to $q_L\oplus q_M$ and $q_S=(q_L\oplus q_M)|(H_S)^\perp/H_S$.  
The primitivity of the embeddings $L\subset S$ and $M\subset S$ is equivalent to the 
condition that the projections $p_L:H_S\to A_L$ and $p_M:H_S\to A_M$ are embeddings.
Let $H_{S,L}=p_L(H_S)$, $H_{S,M}=p_M(H_S)$.  Then $H_S$ is the graph of the isomorphism
$\gamma^L_{L,M}= p_M\circ p_L^{-1}: H_{S,L}\to H_{S,M}$.  
Note that the isotropy of $H_S$ implies that $q_M\circ \gamma^L_{L,M} = -q_L$.   
Conversely, for even lattices $L, M$, subgroups $H_L\subset A_L$, $H_M\subset A_M$, 
and an isomorphism $\gamma: H_L\to H_M$ satisfying $q_M\circ \gamma = - q_L$, define
$S = \{ (x,y) \in L^*\oplus M^* \ : \ x + L \in H_L, \ y+M\in H_M, \ \gamma(x+L) = y+M \}$,
which is an even lattice containing $L, M$ as primitive sublattices. 

Thus, for even lattices $S, L$, the existence of a primitive embedding of $L$ into $S$
is equivalent to  the existence of an even lattice $M$ with $\rank(M)= \rank(S)-\rank(L)$, subgroups $H_L, H_M$ and an isomorphism $\gamma$ satisfying the above conditions
 (\cite[Proposition 1.5.1]{Nikulin2}).
Note that the condition $l(S) \leq l(L)+l(M)$ is necessary. 

To show the existence (or non-existence) of the orthogonal complement $M$ of $L$ in $S$, 
we need to study the discriminant quadratic form $q_M$ of $M$.  
For generators $q_{\theta}^{(p)}(p^k)$, $u_+^{(2)}(2^k)$,
$v_+^{(2)}(2^k)$ of discriminant quadratic forms and their relations, see \cite[Propositions 1.8.1, 1.8.2]{Nikulin2}.  
In the following, we use only the fact that
$q_{D_8} = u_+^{(2)}(2)$, $q_{D_4}=v_+^{(2)}(2)$, $q_{\la \pm 4\ra}=q_{\pm 1}^{(2)}(2^2)$,
where $\la \pm 4\ra$ is the lattice of rank 1 generated by a vector of norm $\pm 4$.
For more details, we refer 
the reader to \cite{Nikulin2}.

\begin{theorem}\label{Artin2}
Let $X$ be a supersingular $K3$ surface with Artin invariant $\sigma$ 
in characteristic two.
Then $X$ contains twenty-eight $(-2)$-curves forming $(12_4, 16_3)$-configuration 
if and only if $\sigma \leq 2$.
\end{theorem}
\begin{proof}
Let $S_\sigma$ be the Picard lattice of a supersingular $K3$ surface with Artin invariant $\sigma$
in characteristic two.  Recall that $S_\sigma$ depends only on $\sigma$ and has signature $(1,21)$, 
$S_\sigma^*/S_\sigma \cong (\bbZ/2\bbZ)^{2\sigma}$, and of type ${\rm I}$ 
in the sense of Rudakov-Shafarevich \cite{RS}, that is, $q_{S_\sigma}(A_{S_\sigma})\subset 
\bbZ/2\bbZ$.
Let $L$ be the lattice generated by twenty-eight $(-2)$-curves which
is isomorphic to $U\oplus E_8\oplus D_8\oplus \la -4\ra$ (see Proposition \ref{PicardLattice}). 
We need to show that $L$ can be primitively embedded into $S_{\sigma}$ if and only if 
$\sigma \leq 2$.

Assume that $L$ can be primitively embedded in $S_\sigma$ and let $M$ be the orthogonal complement of
$L$ in $S_\sigma$, which has rank 3.  Then, $l(L)+l(M) \geq l(S_\sigma)$.  
This implies that $3 + l(M) \geq 2\sigma$, that is, $2\sigma -3 \leq l(M) \leq {\rm rank}(M)=3$
which implies $\sigma \leq 3$.

Conversely, if $\sigma =1$, the assertion follows from the previous theorem.
If $\sigma=2$, 
then $S_\sigma \cong U \oplus E_8\oplus D_8\oplus D_4$.
Obviously, $L=U \oplus E_8\oplus D_8\oplus \la -4\ra$ can be primitively 
embedded in $S_\sigma$.

If $\sigma =3$, then $S_\sigma \cong U \oplus E_8\oplus D_4 \oplus D_4 \oplus D_4$.
In this case, we use a result and the notation of \cite{Nikulin2} to show the non-existence of a primitive embedding of $L$ into $S_\sigma$.  Since $q_{D_4}=v_+^{(2)}(2)$,
$q_{D_8}=u_+^{(2)}(2)$ and $u_+^{(2)}(2)\oplus u_+^{(2)}(2)= v_+^{(2)}(2)\oplus v_+^{(2)}(2)$ 
(\cite[Proposition 1.8.2, b]{Nikulin2}),
$q_{S_\sigma} = u_+^{(2)}(2)\oplus u_+^{(2)}(2)\oplus v_+^{(2)}(2)$.  On the other hand, 
$q_L= u_+^{(2)}(2)\oplus q_{-1}^{(2)}(2^2)$.   Assume that there exists a primitive embedding of $L$ into
$S_{\sigma}$ and let $M$ be the orthogonal complement.  
Since $S_\sigma^*/S_\sigma$ is a 2-elementary abelian group, the discriminant quadratic
form $q_M$ of $M$ is equal to $q = q_{1}^{(2)}(2^2) \oplus v_+^{(2)}(2) (= q_{5}^{(2)}(2^2) \oplus u_+^{(2)}(2))$ because only
the unique isotropic subgroup $\bbZ/2\bbZ$ of $q_{-1}^{(2)}(2^2) \oplus q_{1}^{(2)}(2^2)$ (or 
$q_{-1}^{(2)}(2^2) \oplus q_{5}^{(2)}(2^2)$) can 
define an overlattice $S_\sigma$ of $L\oplus M$.
Note that ${\rm discr}(v_+^{(2)}(2))=12$ and ${\rm discr}(q_{1}^{(2)}(2^2)) =4$ (\cite[Proposition 1.8.1]{Nikulin2}).  
Therefore,
the non-existence of such $M$ follows from \cite[Theorem 1.10.1]{Nikulin2}. 
In fact, the condition (4) of  \cite[Theorem 1.10.1]{Nikulin2} does not hold in our case:
$2^4 = |A_q| (= |M^*/M|) \not\equiv \pm  {\rm discr}(q_M) = \pm 12\cdot 4 =\pm 3\cdot 2^4\ ({\rm mod} 
(\bbZ_2^*)^2)$ because $\pm 3 \notin (\bbZ_2^*)^2$ (\cite[\S 3.3, Remark (1)]{Serre}). Thus, in the case of $\sigma =3$, 
$L$ can not be primitively embedded in $S_\sigma$.
\end{proof}

\begin{remark}\label{K-Schroer}
(suggested by the referee).  First, recall the result in \cite{KS}.  
Let $C$ be a cuspidal curve with the cusp $p$.  
Let $X$ be the minimal resolution of the quotient of the self-product $C\times C$ by
a $\mu_2$-action, which is a supersingular K3 surface in characteristic two.  
The surface $X$ depends on the action of $\mu_2$ leading to the   
two-dimensional family of supersingular K3 surfaces.  
The surface $X$ contains thirty $(-2)$-curves forming the dual graph in 
\cite[Figure 1]{KS}.  Each projection $C\times C\to C$ induces a quasi-elliptic fibration 
$\pi_1, \pi_2$ with five reducible singular fibers of type $\tilde{D}_4$ 
and with four sections.  
The additional fiber of type $\tilde{D}_4$ of $\pi_1$ (resp. $\pi_2$) comes form 
$\{p\}\times C$ (resp. $C\times \{p\}$).  
Thus the Mordell-Weil group of $\pi_i$ is $(\bbZ/2\bbZ)^n$ $(n \geq 2)$.
By applying the Shioda-Tate formula, the discriminant of the Picard lattice 
$\Pic(X)$ is less than or equal to $2^{10}/2^4= 2^6$. This
implies that its Artin invariant $\sigma$ is at most 3.  
In fact, any supersingular K3 surface with Artin invariant 
at most $3$ is isomorphic to $X$ (\cite[Theorem 4.1]{KS}). 

In Remark \ref{from Kummer}, we see that the projection $E\times E\to E$ induces an
elliptic fibration $f_i$ on $\widetilde{\Kum}(E\times E)$ with four singular fibers of type
$\tilde{D}_4$ and {\it eight} sections.  On the other hand,
if $X$ has Artin invariant 3, then the Mordell-Weil group of $\pi_1$ is $(\bbZ/2\bbZ)^2$,
that is, $\pi_1$ has exactly {\it four} sections.  This agrees with Theorem \ref{Artin2}.

Now we need more four sections of $\pi_1$ (and also for $\pi_2$) to get 
a supersingular K3 surface $X$ with Artin invarinant 2.  To do this,  
we specialize the action of $\mu_2$ to the one satisfying the condition 
that the image of the proper transform of the diagonal of $C\times C$ on $X$ is 
such a section (for the equation $(5)$ of the vector field $\delta$ in \cite[p.329]{KS},
put $\lambda_i = \mu_i (i=0,2,4)$).  
By using the translations by the above four sections of $\pi_1$, 
we get the remaining three sections.  
The twenty $(-2)$-curves appeared in four singular fibers of type $\tilde{D}_4$ 
(not the additional one) of $\pi_1$ 
and eight $(=4+4)$ sections form $(12_4, 16_3)$-configuration. 
\end{remark}

The linear system \eqref{lnsystem} defined by the configuration of 12 nodes and 16 lines defines a 
degree two 
map $\varphi$ from $X$ to a quadric surface in $\bbP^3$. The composition with a projection 
$\bbP^1\times \bbP^1 \to \bbP^1$ gives a genus one fibration on $X$.  If $\varphi$ is separable, 
then a general fiber is a double cover of $\bbP^1$ ramified at four points, which is impossible in 
characteristic two.  
Thus, the map $\varphi$ must be inseparable, and the fibration is quasi-elliptic.
This is an analog of the map from the Kummer surface to $\bbP^1\times \bbP^1$.

We conclude that, in characteristic two, there are two types of Kummer surfaces: 
if we keep the elliptic
curve, we get $\Kum(E\times E)$, but lose $(12_4,16_3)$-configuration, and vice versa.
This situation is the same as for Kummer surfaces associated with curves of genus two: if we keep the genus 
two curve, then we 
lose $(16_6)$-configuration and vice versa (\cite{Dolgachev24}, 
\cite{KaKo}).

\end{document}